\documentclass{amsart} 
\usepackage{amsmath,amssymb}

\newcommand{\R}{   \ensuremath{ \mathbb{R} }}
\newcommand{\C}{   \ensuremath{ \mathbb{C} }}
\newcommand{\N}{   \ensuremath{ \mathbb{N} }}

\newcommand{\A}{    \ensuremath{ \mathcal{A}  }} 
\newcommand{\B}{    \ensuremath{ \mathcal{B}  }} 
\newcommand{\F}{    \ensuremath{ \mathcal{F}  }} 
\newcommand{\D}{    \ensuremath{ \mathcal{D}  }} 
\newcommand{\V}{    \ensuremath{ \mathcal{V}  }}

\newcommand{\Aut}{  \ensuremath{ \mathrm{Aut} }}  
\newcommand{\calR}{ \ensuremath{ \mathcal{R}  }}
\newcommand{\calQ}{ \ensuremath{ \mathcal{Q}  }} 
\newcommand{\calP}{ \ensuremath{ \mathcal{P}  }} 
\newcommand{\calO}{ \ensuremath{ \mathcal{O}  }}

\newcommand{\zz}{z,\bar{z}}

\newcommand{\Span}{ \ensuremath{ \mathrm{span}  }} 

\renewcommand{\Im}{\ensuremath{ \mathrm{Im} }}
\renewcommand{\Re}{\ensuremath{ \mathrm{Re} }}

\renewcommand{\bar}[1]{   \ensuremath{ \overline{#1} }}
\renewcommand{\hat}[1]{   \widehat{ #1 }}
\renewcommand{\tilde}[1]{ \widetilde{ #1 }}

\newcommand{\del}{    \partial }
\newcommand{\setof}[2]{     \ensuremath{ \left\{ #1 : #2 \right\} }}

\newcommand{\dd}[2]{    \ensuremath{ \frac{ \displaystyle \del {#1} }{ \displaystyle \del {#2} } }}

\newtheorem{thm}{Theorem}[section]
\newtheorem{lem}[thm]{Lemma}
\newtheorem{prop}[thm]{Proposition}
\newtheorem{cor}[thm]{Corollary}

\theoremstyle{definition}

\newtheorem{example}[thm]{Example}

\theoremstyle{remark}

\begin{document}


\title[Rational jet dependence of formal equivalences]{Rational jet dependence of formal equivalences between real-analytic hypersurfaces in $\C^2$}

\author{R. Travis Kowalski}

\address{Mathematics Department, Colorado College, 14 E.~Cache La Poudre Street, Colorado Springs, CO 80903}

\email{\texttt{tkowalski@coloradocollege.edu}}


\begin{thanks}{2003 {\em Mathematics Subject Classification.} 32H12, 32V20.}
\end{thanks}

\begin{abstract} Let $(M,p)$ and $(\hat{M},\hat{p})$ be the germs of real-analytic $1$-infinite type hypersurfaces in $\C^2$.  We prove that any formal equivalence sending $(M,p)$ into $(\hat{M},\hat{p})$ is formally parametrized (and hence uniquely determined by) its jet at $p$ of a predetermined order depending only on $(M,p)$.  As an application, we use this to examine the local formal transformation groups of such hypersurfaces.
\end{abstract}

\maketitle


\section{Introduction}

A \emph{formal (holomorphic) mapping} $H : (\C^2,p) \to (\C^2,\hat{p})$, with $p, \hat{p} \in \C^2$, is a $\C^2$-valued formal power series
\[
   H(Z) = \hat{p} + \sum_{|\alpha| \ge 1} c_\alpha (Z - p)^\alpha, \quad
   c_\alpha \in \C^2, \quad 
   Z = (Z_1,Z_2).
\] 
The map $H$ is \emph{invertible} if there exists a formal map $H^{-1} : (\C^2,\hat{p}) \to (\C^2,p)$ such that
$H(H^{-1}(Z)) \equiv H^{-1}(H(Z)) \equiv Z$ as formal power series, or equivalently, if the Jacobian of $H$ is nonvanishing at $p$.   We shall denote by $J^k (\C^2,\C^2)_{ p,\hat{p} }$ the \emph{jet space} of order $k$ of (formal) holomorphic mappings $(\C^2,p) \to (\C^2,\hat{p})$, and by $j^k_p (H) \in J^k(\C^2,\C^2)_{ p,\hat{p} }$ the \emph{$k$-jet} of
$H$ at $p$. (See Section \ref{sec:Preliminaries and Basic Definitions} for further details.)

Suppose that $(M,p)$ and $(\hat{M},\hat{p})$ are (germs of) real-analytic hypersurfaces at $p$ and $\hat{p}$ respectively, given by the real-analytic, real-valued local defining functions $\rho(Z,\bar{Z})$ and $\hat{\rho}(Z,\bar{Z})$.  The  formal map $H$ is said to
\emph{take $(M,p)$ into $(\hat{M},\hat{p})$} if
\[
   \hat{\rho} \big( H(Z),\bar{H(Z)} \big) \equiv c(Z,\bar{Z}) \rho(Z,\bar{Z})
\] 
(in the sense of power series) for some formal power series $c(Z,\bar{Z})$; if in addition the formal map is invertible, it is called a \emph{formal equivalence} between $(M,p)$ and $(\hat{M},\hat{p})$, and the germs themselves are called \emph{formally equivalent}.

We wish to study the parametrization and finite determination of invertible formal holomorphic mappings of $\C^2$ taking one real-analytic hypersurface $M$ into another.  There is a great deal of literature on this if $M$ is assumed to be
\emph{minimal} at $p$, i.e.~if there is no complex hypersurface through $p$ in $\C^2$ contained in $M$; see the remarks at the end of this introduction.  In the present paper, however, we shall assume that $M$ is not minimal at $p$, so that there exists a complex hypersurface $\Sigma \subset \C^2$ with $p \in \Sigma \subset M$.  It is well known (see \cite{ChernMoser:RealHypersurfaces}; also \cite{BaouendiEbenfeltRothschild:RealSubmanifolds}, Chapter IV) that for any
real-analytic hypersurface $M \subset \C^2$ and point $p \in M$ (not necessarily minimal), there exist local holomorphic coordinates $(z,w) \in \C \times \C$, vanishing at $p$, such that $M$ is defined locally by the equation
\[
    \Im \, w = \Theta( \zz, \Re \, w ),
\] 
where $\Theta(\zz,s)$ is a real-valued, real-analytic function such that
\[
   \Theta(z,0,s) \equiv \Theta(0,\bar{z},s) \equiv 0.  
\] 
Such coordinates are called \emph{normal coordinates} for $M$ at $p$, and are not unique.  $M$ is said to be of \emph{finite type} at $p$ if $\Theta(\zz,0) \not \equiv 0$; otherwise $M$ is of \emph{infinite type} at $p$.  This definition is equivalent to being of finite type in the sense of Kohn \cite{Kohn:BoundaryBehavior} and Bloom and Graham \cite{BloomGraham:OnTypeConditions}.  For real-analytic hypersurfaces, it is also equivalent to minimality --- indeed, if $M$ is of infinite type at $p$, then (in normal coordinates) $M$ contains the nontrivial complex hypersurface $\Sigma = \{ w = 0 \}$.   (See e.g.~\cite{BaouendiEbenfeltRothschild:RealSubmanifolds}, Chapter I, for further details.)

In this paper, we shall focus our attention on $1$-infinite type points $p$ of a real-analytic hypersurface $M \subset \C^2$, i.e.~points at which the normal coordinates above satisfy the additional condition that $\Theta_s(\zz,0) \not\equiv 0$. (See Section \ref{sec:Preliminaries and Basic Definitions} for precise definitions.)  Our main result gives rational dependence of a formal equivalence between $1$-infinite type hypersurfaces on  its jet of a predetermined order.

\begin{thm} \label{THEOREM 1}
Let $M \subset \C^2$ be a real-analytic hypersurface, and suppose $p \in M$ is of $1$-infinite type.  Then there exists an integer $k$ such that given any hypersurface $\hat{M} \subset \C^2$ with $(\hat{M},\hat{p})$ formally equivalent to $(M,p)$, there exists a formal power series of the form
\begin{equation} \label{eqn:form of Psi} 
   \Psi(Z; \Lambda) = \sum_{\alpha} \frac{ p_\alpha (\Lambda) }{ q(\Lambda)^{\ell_\alpha} } (Z - p)^\alpha,
\end{equation} 
where $p_\alpha, q$ are polynomials on the jet space $J^k(\C^2,\C^2)_{ p,\hat{p} }$ valued in $\C^2$ and $\C$ respectively, and the $\ell_\alpha$ are nonnegative integers, such that for any formal equivalence $H : (M,p) \to (\hat{M},\hat{p})$, the following holds:
\[
   q \big( j^k_p (H)  \big) =  \det \left( \dd{H}{Z} (p) \right)  \neq 0, 
   \qquad  \mbox{and}  \qquad
   H(Z) = \Psi \big( Z; j^k_p (H)  \, \big).
\]
\end{thm}

\noindent Our proof (presented in Section \ref{Proofs of the Main Results}) will actually give a constructive process for determining such an $k$.

Theorem \ref{THEOREM 1} has a number of applications.  The first states that any formal equivalence between two germs of $1$-infinite type hypersurfaces $(M,p)$ and $(\hat{M},\hat{p})$ is determined by finitely many derivatives at $p$.

\begin{thm} \label{THEOREM 2}
Let $(M,p)$ and $k$ be as in Theorem \ref{THEOREM 1}.  If $H^1, H^2 : (M,p) \to (\hat{M},\hat{p})$ are formal equivalences and
\[
   \dd{^{|\alpha|} H^1}{Z^\alpha}(p) = \dd{^{|\alpha|} H^2}{Z^\alpha}(p), \quad \forall \, |\alpha| \le k,
\] 
then $H^1 = H^2$ as power series.
\end{thm}

Our second application deals with the structure of jets of formal equivalences in the jet space $J^k(\C^2,\C^2)_{ p, \hat{p} }$, or rather in the submanifold $G^k(\C^2)_{ p,\hat{p} }$ of jets of invertible maps taking $(\C^2,p)$ to $(\C^2,\hat{p})$.  We shall denote by ${\F}(M,p;\hat{M},\hat{p})$ the set of formal equivalences taking $(M,p)$ into $(\hat{M},\hat{p})$.  We have the following.

\begin{thm} \label{THEOREM 3}
Let $(M,p)$ and $k$ be as in Theorem \ref{THEOREM 1}.  Then for any (germ of a) real-analytic hypersurface
$(\hat{M},\hat{p})$ in $\C^2$, the mapping
\[
   j^k_p : {\F}(M,p;\hat{M},\hat{p}) \to G^k(\C^2)_{ p,\hat{p}}
\]
is an injection onto a real algebraic submanifold of $G^k(\C^2)_{ p,\hat{p}}$.
\end{thm}

Of special interest is the case $(\hat{M},\hat{p}) = (M,p)$, since ${\F}(M,p;\hat{M},\hat{p})$ becomes a group under composition, called the \emph{formal stability group} of $M$ at $p$ and denoted by ${\Aut}(M,p)$.  We shall denote by 
$G^k(\C^2)_{p} := G^k(\C^2)_{p,p}$ the \emph{$k$-jet group} of $\C^2$ at $p$.  The following result is then a corollary of Theorem \ref{THEOREM 3}.

\begin{thm} \label{THEOREM 4}
Let $(M,p)$ and $k$ be as in Theorem \ref{THEOREM 1}.  Then the mapping
\[
   j^k_p : {\Aut}(M,p) \to G^k(\C^2)_{p}
\]
defines an injective group homomorphism onto a real algebraic Lie subgroup of $G^k(\C^2)_{p}$.
\end{thm}

The study of the (formal) transformation groups of hypersurfaces in $\C^N$ has a long history.  Its roots can be traced back to E.~Cartan, who studied the structure of the local transformation groups of smooth Levi nondegenerate hypersurfaces in $\C^2$ in \cite{ECartan:SurLaGeometrieI}, \cite{ECartan:SurLaGeometrieII}.  These results were later extended to higher dimensions by Chern and Moser in \cite{ChernMoser:RealHypersurfaces}, who also proved the finite determination of such equivalences by their $2$-jets.

Further results about the  transformation groups of various classes of finite type generic submanifolds of $\C^N$ were more recently  obtained by a number of mathematicians.  We mention here the work of Tumanov and Henkin \cite{TumanovHenkin:LocalCharacterization}, Tumanov \cite{Tumanov:FiniteDimensionality}, Zaitsev \cite{Zaitsev:GermsOfLocalAutomorphisms}, and Baouendi, Ebenfelt, and Rothschild
\cite{BaouendiEbenfeltRothschild:Parametrization}, \cite{BaouendiEbenfeltRothschild:RationalDependence}. In
\cite{BaouendiEbenfeltRothschild:Parametrization}, the authors showed that modified versions of Theorems \ref{THEOREM 2}--\ref{THEOREM 4} hold for smooth hypersurfaces $M,\hat{M}$ in $\C^N$ with $M$ of finite type and $\hat{M}$ finitely  nondegenerate.  This was later extended to smooth generic submanifolds in \cite{Zaitsev:GermsOfLocalAutomorphisms} and \cite{BaouendiEbenfeltRothschild:RationalDependence}.  In particular, in the latter paper, the authors proved an analogue of Theorem \ref{THEOREM 1} for generic submanifolds $M$ and $\hat{M}$ of finite type and finitely nondegenerate.  They
also proved the convergence of the formal mappings as well.

For the proofs of the four theorems above, it is convenient to work with formal mappings between formal real hypersurfaces.  Hence, the results presented here will be reformulated  and proved  in this more general context.  The following section presents the necessary preliminaries and definitions.  In what follows, the distinguished points $p$ and $\hat{p}$ on $M$ and $\hat{M}$, respectively, will, for convenience and without loss of generality, be assumed to be $0$.


\section{Preliminaries and basic definitions}
\label{sec:Preliminaries and Basic Definitions}

\subsection{Formal mappings and hypersurfaces}
\label{sec:Formal Mappings and Manifolds}

Let $X = (X_1,\dots,X_N)$ denote a $N$-tuple of indeterminates, and let $\calR$ denote a commutative ring with unity.  We shall define the following rings:
\begin{itemize}
   \item $\calR[[X]] :=$ the ring of formal power series in $X$ with coefficients in $\calR$.
   \item $\calR[X] :=$ the ring of polynomials in $X$ with coefficients $\calR$.
\end{itemize} 
For $\calR = \C$, we shall also define the rings
\begin{itemize}
   \item $\C\{X\} :=$ the ring of convergent power series in $X$ with coefficients in $\C$.
   \item $\calO_\epsilon(X) :=$ the ring of power series in $X$ with  coefficients in $\C$ which converge for
$X_j \in \C$, $|X_j| < \epsilon$, $1 \le j \le N$.
\end{itemize} 
Observe that we have canonical embeddings
\[
   \C[X] \subset \calO_\epsilon(X) \subset \C\{X\} \subset \C[[X]].
\]

A power series $\rho \in \C[[Z,\zeta]]$, where $Z = (Z_1,\dots,Z_N)$ and $\zeta = (\zeta_1,\dots,\zeta_N)$, is called \emph{real} if $\rho(Z,\zeta) = \bar{\rho}(\zeta, Z)$, where $\bar{\rho}$ denotes the power series obtained by replacing the
coefficients of $\rho$ by their complex conjugates.  If, in addition, the power series $\rho$ satisfies the two conditions
\begin{equation} \label{eqn:wedge}
   \rho(0)=0, \quad d \rho (0) \neq 0,
\end{equation}
then we say $\rho$ defines a \emph{formal real hypersurface} $M$ of $\C^N$ through $0$, and we shall symbolically write
\[
   M = \big\{ \rho \big( Z, \bar{Z} \, \big) = 0 \big\},
\]
and say that the pair $(M,0)$ is a \emph{formal real hypersurface}.  The function $\rho$ is a \emph{formal defining function} for $M$.  The reader should observe that if $M$ is a formal real hypersurface in $\C^N$ with formal defining function $\rho$, then in general there is no actual point set $M \subset \C^N$.

Suppose $\hat{\rho}$ is another formal power series (not necessarily real) which satisfies the conditions (\ref{eqn:wedge}).  If there exists a power series $a(Z,\zeta)$ (necessarily invertible at $0$) such that
\[
   \hat{\rho}(Z,\zeta) = a(Z,\zeta) \, \rho(Z,\zeta),
\] 
then we say that $\hat{\rho}$ also defines the formal real hypersurface $M$, and shall also write $M = \{ \hat{\rho}(Z,\bar{Z})=0 \}$.

By a formal mapping $H : (\C^N,0) \to (\C^N,0)$, denoted $H \in {E}(\C^N,\C^N)_{0,0}$, we shall mean an element $H \in \C[[Z]]^{N}$ such that $H(0) = 0$.  We say $H$ is a \emph{formal change of coordinates} if it is formally invertible, i.e.~if there exists a formal map $H^{-1} : (\C^N,0) \to (\C^N,0)$ such that
\[
   H(H^{-1}(Z)) \equiv H^{-1}(H(Z)) \equiv Z
\]
as formal power series. As noted in the introduction, $H$ is a formal change of coordinates in $\C^N$ if and only if its Jacobian at $0$ is nonzero.

Given a formal change of coordinates $H$ in $\C^N$, we define its corresponding \emph{formal holomorphic change of
variable} by
\[
   Z = H(Z'), \quad \zeta = \bar{H}(\zeta').
\]
If $M = \{ \rho(Z,\bar{Z})=0 \}$ is a formal real hypersurface of $\C^N$, then we say $M$ is expressed in the $Z'$
coordinates by $\{ \rho(H(Z'),\rho(\bar{H(Z')}) = 0 \}$.

If $\hat{M} = \{ \hat{\rho}(Z,\bar{Z}) = 0\}$ is another formal real hypersurface of $\C^{N}$, then a formal mapping $H \in {E}(\C^N,\C^N)_{0,0}$ is said to \emph{take $M$ into $\hat{M}$}, denoted $H : (M,0) \to (\hat{M},0)$, if there exists a 
power series $c(Z,\zeta)$ (not necessarily invertible at $0$) such that
\[
   \hat{\rho} \big( H(Z), \bar{H}(\zeta) \big)
   = c(Z,\zeta) \, \rho(Z,\zeta).
\] 
This definition is independent of the power series used to define $M$ and $\hat{M}$.

If $H : (M,0) \to (\hat{M},0)$ is as above, and $H$ is invertible, then it follows that $H^{-1}$ takes $\hat{M}$ into $M$.  In this case, we say that $M$ and $\hat{M}$ are \emph{formally equivalent}, and that $H$ is a \emph{formal equivalence} between them, denoted $H \in {\F}(M,0;\hat{M},0)$.

The motivation behind these definitions is the following.  If the formal series $\rho$ defining the formal real hypersurface $M$ is actually convergent, then the equation $\rho(Z,\bar{Z}) = 0$ defines a real-analytic hypersurface $M$ of $\C^N$ passing through the origin.  Moreover, if $H : \C^N \to \C^N$ is a holomorphic mapping with $H(0) = 0$, and $M,\hat{M}$ are both real-analytic hypersurfaces of $\C^N$, then $H(M) \subset \hat{M}$ if and only if the formal mapping $H$ maps the formal real hypersurface $M$ into the formal real hypersurface $\hat{M}$.

For each positive integer $k$, we denote by $J^k(\C^N,\C^N)_{0,0}$ the jet space of order $k$ of (formal) holomorphic mappings $(\C^N,0) \to (\C^N,0)$, and by $j^k_0 : {E}(\C^N,\C^N) \to J^k(\C^N,\C^N)_{0,0}$ the corresponding jet mapping taking a formal mapping $H$ to its $k$-jet at $0$, $j^k_0 (H)$.  We shall denote by $G^k(\C^N)_0 \subset J^k(\C^N,\C^N)_{0,0}$ the collection of $k$-jets of invertible formal mappings of $(\C^N,0)$ to itself.

Given coordinates $Z$ and $\hat{Z}$ on $\C^N$, we may identify the jet space $J^k(\C^N,\C^N)_{0,0}$ with the set of degree $k$ polynomial mappings of $(\C^N,0) \to (\C^N,0)$.  The coordinates on $J^k(\C^N,\C^N)_{0,0}$, which we shall denote by $\Lambda$, can then be taken to be the coefficients of these polynomials.  Observe that formal changes of coordinates in $\C^N$ yield polynomial changes of coordinates in $J^k(\C^N,\C^N)_{0,0}$.

If $M$ is a formal real hypersurface in $\C^N$, then there is a formal change of coordinates $Z = (z,w) \in \C[[z,w]]^{N}$ with $z=(z_1,\dots,z_{N-1})$, such that $M$, under the corresponding formal holomorphic change of variable $Z = Z(z,w)$, $\zeta = \bar{Z}(\chi,\tau)\big)$, is defined by
\[
   \rho(z,w,\chi,\tau) := 
   \bigg( \frac{w - \tau}{2 i} \bigg) - \Theta \bigg( z, \chi, \frac{w + \tau}{2} \bigg)
   \in \C[[Z,\zeta]],
\]
where $\Theta \in \C[[z,\chi,s]]$ is real and satisfies $\Theta(z,0,s) = \Theta(0,\chi,s) = 0$.  Such coordinates are called \emph{normal coordinates} for $M$. (See \cite{BaouendiEbenfeltRothschild:RealSubmanifolds}, Chapter IV.)

Using the formal Implicit Function Theorem to solve for $w$ above, there exists a unique formal power series $Q \in \C[[z,\chi,\tau]]$ with $Q(0,0,0) = 0$ such that $\rho \big( z, Q(z,\chi,\tau), \chi, \tau \big) \equiv 0$; moreover, $Q$ is convergent whenever $\Theta$ is.  In particular, this implies that there exists a  power series $a(z,w,\chi,\tau)$,  nonvanishing at $0$, such that
\[
   \rho(z,w,\chi,\tau) = a(z,w,\chi,\tau) \big( w - Q(z,\chi,\tau) \big),
\]
whence we may write (under an abuse of notation)
\begin{equation} \label{eqn:I in Q}
   M = \bigg\{ \bigg( \frac{w - \bar{w}}{2 i} \bigg) = 
       \Theta \bigg( z, \bar{z}, \frac{w + \bar{w}}{2} \bigg) \bigg\} 
     = \big\{ w = Q(z,\bar{z},\bar{w}) \big\}.
\end{equation}
Observe that the normality of the coordinates implies $Q(z,0,\tau) = Q(0,\chi,\tau) = \tau$.

Given normal coordinates $Z = (z,w)$ for $M$ as above, define the numbers $m,r,L,K \in \{0,1,2,\cdots\} \cup \{
\infty \}$ as follows.  Set
\begin{equation} \label{eqn:m defined}
   m   := \sup \setof{ q }{ \Theta_{s^j}(z,\chi,0) \equiv 0 \quad \forall \, j < q }
\end{equation}
If $m = \infty$ (i.e.~if $\Theta \equiv 0$), then set $r = L = K = \infty$.  Otherwise, set
\begin{align} 
\label{eqn:r defined}
   r &:= \sup \setof{ q }{ \Theta_{z^\alpha \chi^\beta s^m}(0,0,0) = 0 
         \quad \forall \, |\alpha| + |\beta| < q }, \\
\label{eqn:L defined}
   L &:= \sup \setof{ q }{ \Theta_{\chi^\beta s^m}(z,0,0) \equiv 0 
         \quad \forall \, |\beta| < q }, \\
\label{eqn:K defined}
   K &:= \sup \setof{ q }{ \Theta_{z^\alpha \chi^\beta s^m}(0,0,0) = 0 
         \quad \forall \, |\alpha| < q, \, |\beta| = L}.
\end{align}
We shall show (Theorem \ref{4-tuple}) that this $4$-tuple of numbers is independent of the normal coordinates used to define them.

We say that $M$ is of \emph{finite type} at $0$ if $m = 0$; otherwise  $M$ is  of \emph{infinite type} at $0$.  If we wish to emphasize the number $m \ge 1$, we shall say that $M$ is of \emph{$m$-infinite type} at $0$ if $m < \infty$, and is \emph{flat} at $0$ if $m = \infty$.  We shall further say $M$ is of \emph{finite type $r$} at $0$ if $m = 0$, and is of \emph{$m$-infinite type $r$} at $0$ if $1 \le m < \infty$.

We conclude these definitions by stating a few known results concerning these numbers in the case when $M$ is a real-analytic hypersurface in $\C^N$.  In this case, it is known that the pair $(m,r)$ is a biholomorphic invariant of $M$; see \cite{Meylan:AReflectionPrinciple}.   If $m > 0$, i.e.~$M$ is of infinite type at $0$, then $M$ contains a formal complex hypersurface $\Sigma$ passing through $0$.  (Indeed, in normal coordinates, we may take $\Sigma = \{ w = 0 \}$.)  It is also known that if $m > 0$ at the origin, then $m$ is actually constant along the complex hypersurface $\Sigma \subset M$ through $0$.  And while $r$ is not constant along $\Sigma$, it is known that there exists a proper, real-analytic subvariety $V \subset \Sigma$ outside of which all points are of $m$-infinite type $2$.  We direct the reader to \cite{Ebenfelt:AnalyticityNonminimal} for more details.

\subsection{Statement of results}
\label{sec:Statement of Results}

Our first result shows that the $4$-tuple $(m,r,L,K)$ (and hence the notion of being $m$-infinite type $r$ at a point) is in fact a \emph{formal} invariant of a hypersurface.  Specifically, we have the following.

\begin{thm} \label{4-tuple}
Let $(M,0)$ be a formal real hypersurface of $\C^{N}$. Then the numbers $(m,r,L,K)$ are independent of the choice of normal coordinates used to define them.  Moreover, if $(\hat{M},0)$ is formally equivalent to $(M,0)$ and has the corresponding $4$-tuple
$(\hat{m},\hat{r},\hat{L},\hat{K})$, then $(m,r,L,K) = (\hat{m},\hat{r},\hat{L},\hat{K})$.
\end{thm}

We shall then focus exclusively on the case $N=2$ and $m=1$.  We may now state the generalizations of Theorems \ref{THEOREM 1} through \ref{THEOREM 4} valid for formal real hypersurfaces.  Our main result is the following.

\begin{thm} \label{Parametrization}
Let $(M,0)$ be a formal real hypersurface in $\C^2$ which is of $1$-infinite type.  Then there exists an integer $k$ such that given any formal real hypersurface $(\hat{M},0)$ in $\C^2$ formally equivalent to $(M,0)$, there exists a formal power series of the form
\begin{equation} \label{eqn:form of Psi II}
   \Psi(Z; \Lambda) = \sum_{\alpha} \frac{ p_\alpha (\Lambda) }{ q(\Lambda)^{\ell_\alpha} } Z^\alpha,
\end{equation}
where $p_\alpha, q$ are polynomials on the jet space $J^k(\C^2,\C^2)_{0,0}$ valued in $\C^2$ and $\C$ respectively,  and the $\ell_\alpha$ are nonnegative integers, such that for any formal equivalence $H \in {\F}(M,0;\hat{M},0)$, the following holds: 
\[
   q \big( j^k_0 (H) \big) =  \det \left( \dd{H}{Z} (0) \right) \neq 0, \qquad
   H(Z) = \Psi \big( Z;  j^k_0 (H)  \, \big).
\]
\end{thm}

It is clear from the remarks made in the previous section that Theorem \ref{Parametrization} is a more general version of Theorem \ref{THEOREM 1} from the introduction.  As a consequence of this result, we have the following, from which Theorem \ref{THEOREM 2} is derived.

\begin{thm} \label{Finite Determination}
Let $(M,0)$ be a formal real hypersurface in $\C^2$ of $1$-infinite type, and let $k$ be the number described in Theorem \ref{Parametrization}.  Then for any  formal hypersurface $(\hat{M},0)$ formally equivalent to $(M,0)$, and any formal equivalences $H^1, H^2 : (M,0) \to (\hat{M},0)$, if
\[
   \dd{^{|\alpha|} H^1}{Z^\alpha}( 0 ) = \dd{^{|\alpha|} H^2}{Z^\alpha}( 0 ) \quad \forall \, |\alpha| \le k,
\]
then $H^1 = H^2$ as power series.
\end{thm}

We shall then prove the following generalization of Theorem \ref{THEOREM 4}.

\begin{thm} \label{Lie Subgroup} 
Let $M$ and $k$ be as in Theorem \ref{Parametrization}.  Then the mapping 
\[
   j^k_0 : {\Aut}(M,0) \to G^k(\C^2)_{0}
\]
defines an injective group homomorphism onto a real algebraic Lie subgroup of $G^k(\C^2)_{0}$.
\end{thm}

A consequence of Theorem \ref{Lie Subgroup} is the following, which is a generalization of Theorem \ref{THEOREM 3}.

\begin{thm} \label{Varieties} Let $M$ and $k$ be as in Theorem \ref{Parametrization}.  Then for any formal real hypersurface $\hat{M}$ in $\C^2$, the mapping
\[
   j^k_0 : {\F}(M,0;\hat{M},0) \to J^k(\C^2)_{0}
\]
is an injection onto a real algebraic submanifold of $G^k(\C^2)_{0}$.
\end{thm}


\section{Formal invariance of type conditions}
\label{sec:Formal Invariance of Type Conditions}

In this section, we shall prove Theorem \ref{4-tuple}.  In fact, we shall prove the following, slightly sharper statement, of which Theorem \ref{4-tuple}  is an immediate consequence.

\begin{prop} \label{invariant pair}
Let $(M,0)$ be a formal real hypersurface in $\C^{N}$, given in normal coordinates $Z=(z,w)$ by equation (\ref{eqn:I in Q}).  Let $(\hat{M},0)$ be a formal real hypersurface in $\C^{N}$, given in normal coordinates $\hat{Z}=(\hat{z},\hat{w})$ by the corresponding ``hatted'' defining functions, i.e.~of the form
\[
   \hat{M} =  \bigg\{ \frac{\hat{w} - \bar{ \hat{w} }}{2 i} 
     = \hat{\Theta} \bigg( \hat{z}, \bar{\hat{z}}, \frac{\hat{w} + \bar{ \hat{w} }}{2} \bigg) \bigg\}
     = \bigg\{ \hat{w} = \hat{Q} \big( \hat{z}, \bar{\hat{z}}, \bar{\hat{w}} \big) \bigg\}.
\]
Define the $4$-tuple $(m,r,L,K)$ for $M$ as in Section \ref{sec:Preliminaries and Basic Definitions}, and define the corresponding $4$-tuple $(\hat{m},\hat{r},\hat{L},\hat{K})$ for $\hat{M}$.  If $M$ and $\hat{M}$ are formally equivalent, then $(m,r,L,K) = (\hat{m},\hat{r},\hat{L},\hat{K})$.
\end{prop}

We begin with a useful lemma concerning the form of formal mappings in normal coordinates.  It is proved as Lemma 9.4.4 in \cite{BaouendiEbenfeltRothschild:RealSubmanifolds}, Chapter IX.

\begin{lem} \label{form of H}
Let $M, \hat{M}$ be formal hypersurfaces in $\C^N$ through $0$, expressed in normal coordinates as in Proposition \ref{invariant pair}.  If $H = (F,G): (M,0) \to (\hat{M},0)$ is a formal mapping, then  $G(z,w) = w \, g(z,w)$ for some $g \in \C[[z,w]]$.  Moreover, if $H$ is a formal equivalence, then $F(z,0) \in \C[[z]]^{N-1}$ is a formal equivalence, and $g(0,0) \neq 0$.
\end{lem}

As a consequence of this lemma, we shall henceforth write formal equivalences (in suitable normal coordinates) as
\begin{equation} \label{eqn:H expanded}
   H(z,w) = \big( f(z,w), w \, g(z,w) \big),
\end{equation}
with $f = (f^1, \dots, f^{N-1} ) \in \C[[z,w]]^{N-1}$ satisfying $\det f_z(0,0) \neq 0$ and  $g \in \C[[z,w]]$ satisfying $g(0,0) \neq 0$.  Observe that the condition that $H$ map $M$ formally into $\hat{M}$ may be written as
\begin{equation} \label{eqn:formal map, Q}
   Q(z,\chi,\tau) \, g \big( z, Q(z,\chi,\tau) \big)
   \equiv \hat{Q} \big( f \big( z, Q(z,\chi,\tau) \big), \bar{f}(z,\chi), \tau \, \bar{g}(\chi,\tau) \big).
\end{equation}
Moreover, for convenience, we shall formally expand $f$ and $g$ as
\begin{equation} \label{eqn:fn and gn expanded}
   f(z,w) = \sum_{n\ge 0} \frac{f_n(z)}{n!} \, w^n, \quad
   g(z,w) = \sum_{n\ge 0} \frac{g_n(z)}{n!} \, w^n.
\end{equation}

The main technical lemma in the proof of Proposition \ref{invariant pair} is the following.

\begin{lem} \label{m' to m}
Suppose that $M, \hat{M}$ are formal hypersurfaces in $\C^N$ through $0$, expressed in normal coordinates as in Proposition \ref{invariant pair}, and assume that $H : (M,0) \to (\hat{M},0)$ is a formal equivalence.  Then for every $j \ge 0$, if
\[
   \hat{Q}(\hat{z},\hat{\chi},0) \equiv \hat{Q}_{\hat{\tau}}(\hat{z},\hat{\chi},0) - 1 
      \equiv \hat{Q}_{{\hat{\tau}}^2}(\hat{z},\hat{\chi},0) \equiv
      \cdots \equiv \hat{Q}_{{\hat{\tau}}^j}(\hat{z},\hat{\chi},0) \equiv 0,
\]
then
\begin{equation} \label{eqn:lotsa tilde Q's}
   Q(z,\chi,0) \equiv Q_{\tau}(z,\chi,0) - 1 \equiv Q_{{\tau}^2}(z,\chi,0) \equiv
      \cdots \equiv Q_{{\tau}^j}(z,\chi,0) \equiv 0.
\end{equation}
Moreover,  $g_0(z), g_1(z), \dots, g_j(z)$ are all real constants (with $g_0(z)$  nonzero), and
\[
   Q_{\tau^{j+1}}(z, \chi,0) \equiv g(0)^j \, \hat{Q}_{{\hat{\tau}}^{j+1}} \big( f_0(z), \bar{f_0}(\chi), 0 \big).
\]
\end{lem}

To prove Lemma \ref{m' to m}, we shall make use of the following two results.  The first result is a generalization of the Chain Rule due to Faa de Bruno; see e.g.~\cite{Range:HolomorphicFunctions}:

\begin{lem}[Faa de Bruno's Formula] \label{Faa de Bruno}
Suppose that $f = \big( f_1, f_2, \dots, f_\ell \big) \in \C^\ell[[z]]$ with $z \in \C$ and $f(0) = 0$, and suppose $h(z_1,z_2, \dots, z_\ell) \in \C[[z_1,z_2,\dots,z_\ell]]$.  Then
\begin{align*}
   \dd{^v}{z^v} \big\{ h \big( f(z) \big) \big\}
   &= \sum_{ \substack{
      [\alpha^1]+[\alpha^2]+\cdots  \\ 
      + [\alpha^\ell] = v } } 
        \frac{v! \, h_{{z_1}^{|\alpha^1|} {z_2}^{|\alpha^2|} \cdots {z_\ell}^{|\alpha^\ell|}} \big( f(z) \big)}{
              \alpha^1! \, \alpha^2! \cdots \alpha^\ell!} 
        \prod_{ \substack{ 
           1 \le q \le v \\ 
           1 \le p \le \ell } } 
          \bigg( \frac{{f_p}^{(q)}(z)}{q!} \bigg)^{\alpha^p_q},
\end{align*}
where each $\alpha^p = (\alpha^p_1, \dots, \alpha^p_v)$ denotes an $v$-dimensional multi-index, and
\[
   |\alpha^p| = \sum_{q=1}^v \alpha^p_q, \quad
   [\alpha^p] = \sum_{q=1}^v q \, \alpha^p_q, \quad
    \alpha^p! = \prod_{q=1}^v (\alpha^p_q)!
\]
\end{lem}

The proof is a routine induction, and is left to the reader.  The other result we shall need gives a second characterization of the number $m$.  It is proved as Proposition 1.7 in \cite{BaouendiEbenfeltRothschild:ReflectionPrinciple}.

\begin{prop} \label{m, equivalently}
Let $M$, $m$, $\Theta$, and $Q$ be as above.  Then
\[
   m = \sup \setof{ q }{ \dd{^j}{\tau^j} \big\{ Q(z, \chi, \tau) - \tau \big\} \bigg|_{\tau = 0} \equiv 0
       \quad \forall \, j < q }.
\]
Furthermore,
\[
   Q_{\tau^m}(z,\chi,0) = \left\{ \begin{array}{l @{\quad} l}
      \frac{ \displaystyle 1 + i \, \Theta_s(z,\chi,0)}{\displaystyle 1 - i \, \Theta_s(z,\chi,0)} & m = 1 \\
           & \\
      2 i \, \Theta_{s^m}(z,\chi,0) & 2 \le m < \infty
   \end{array} \right. .
\]
\end{prop}

We now proceed with the proof of Lemma \ref{m' to m}.

\begin{proof}
To begin, observe that differentiating identity (\ref{eqn:formal map, Q}) $v$ times in $\tau$, setting $\tau = 0$, and canceling a $v!$ from both sides yields the identity
\begin{align} \label{eqn:formal map, Q tau}
   &\sum_{k + [\xi]=v} \frac{ g_{|\xi|}(z) \, Q_{\tau^k}(z,\chi,0) }{k! \, \xi!} \prod_{p=1}^v \bigg( \frac{Q_{\tau^p}(z,\chi,0)}{p!} \bigg)^{\xi_p}  \\
   &\quad \equiv \sum_{ \substack{
       [\alpha^1]+\cdots+[\alpha^n]+[\beta^1]+\cdots \\
       \cdots+[\beta^n]+[\gamma]=v } }
   \! \! \! \frac{\hat{Q}_{\hat{z}^{(|\alpha^1|,\dots,|\alpha^n|)} \hat{\chi}^{(|\beta^1| \cdots, |\beta^n|)} \hat{\tau}^{|\gamma|}} \big( f_0(z), \bar{f_0}(\chi), 0 \big) }{
   \alpha^1!\cdots \alpha^n! \, \beta^1! \cdots \beta^n! \, \gamma!}  \nonumber \\
  &\, \, \times \! \! \prod_{ \substack{ 
      1\le q \le v \\ 
      1\le u \le n } } 
        \! \left( \sum_{[\eta]=q} \frac{f^u_{|\eta|}(z)}{\eta!}  \prod_{r=1}^q \bigg( \frac{Q_{\tau^r}(z,\chi,0)}{r!} \bigg)^{\eta_r} \right)^{\alpha^u_q}
        \! \bigg( \frac{\bar{f^u_q}(\chi)}{q!} \bigg)^{\beta^u_q}
        \! \bigg( \frac{\bar{g_{q-1}}(\chi)}{(q-1)!} \bigg)^{\gamma_q}. \nonumber
\end{align}

We now proceed by induction.  For $j = 0$, we assume only that $\hat{Q}(\hat{z},\hat{\chi},0) \equiv 0$.  Setting $\tau = 0$ in identity (\ref{eqn:formal map, Q}), we find
\[
   Q(z,\chi,0) \ g \big( z, Q(z,\chi,0) \big)
   \equiv \hat{Q} \big( f_0(z), \bar{f_0}(\chi),0 \big) = 0.
\]
Since $g(z,Q(z,\chi,0))$ does not vanish at $z = \chi = 0$, we conclude $Q(z,\chi,0) \equiv 0$.

Applying the $v=1$ case of identity (\ref{eqn:formal map, Q tau}), we find
\[
   Q_\tau(z,\chi,0) g_0(z) \equiv \hat{Q}_{\hat{\tau}} \big( f_0(z), \bar{f_0}(\chi),0 \big) \bar{g_0}(\chi).
\]
Setting $\chi = 0$ yields $g_0(z) \equiv \bar{g_0}(0) = \bar{g_0(0)}$, whence $g_0(z)$ is a real constant $r$, and since $H$ is invertible, $r\neq 0$ necessarily.  Dividing $g_n(z) = \bar{g_0}(\chi) = r \neq 0$ from both sides of the identity above yields
\[
   Q_\tau(z,\chi,0) \equiv \hat{Q}_{\hat{\tau}} \big( f_0(z), \bar{f_0}(\chi),0 \big),
\]
which proves the $j = 0$ case.

Now, assume that the lemma holds for some $j-1 \ge 0$; we shall prove it for $j$.  Suppose that (\ref{eqn:lotsa tilde Q's}) holds.  By induction, we know that
\[
   Q(z,\chi,0) \equiv Q_{\tau}(z,\chi,0) - 1 \equiv Q_{{\tau}^2}(z,\chi,0) \equiv
   \cdots \equiv Q_{{\tau}^{j-1}}(z,\chi,0) \equiv 0,
\]
that $g_0, g_1, \dots, g_{j-1}$ are constant functions, and that
\[
   Q_{\tau^j}(z,\chi,0) \equiv r^{j-1} \, \hat{Q}_{{\hat{\tau}}^j} \big( f_0(z), \bar{f_0}(\chi), 0 \big).
\]
In the $j = 1$ case, this implies $Q_\tau(z,\chi,0) \equiv 1$; otherwise it implies $Q_{\tau^j}(z,\chi,0) \equiv 0$, as desired.

Substituting these values into identity (\ref{eqn:formal map, Q tau}) (with $v = j+1$), we obtain
\[
   r \, Q_{\tau^{j+1}}(z,\chi,0) + (j+1) g_j(z)
   \equiv r^{j+1} \hat{Q}_{{\hat{\tau}}^{j+1}} \big( f_0(z), \bar{f_0}(\chi), 0 \big)  + (j+1) \bar{g_j}(\chi).
\]
Setting $\chi = 0$ yields
\[
   (j+1) g_j(z) = (j+1) \bar{g_j}(0) = (j+1) \bar{g_j(0)},
\]
so $g_j(z)$ is a real constant.  Subtracting $(j+1) g_j(z)$ from both sides and dividing by $r\neq 0$ completes the induction.
\end{proof}

\begin{cor} \label{m as invariant}
Let $M,\hat{M}$ be formal real submanifolds of $\C^N$ through $0$, given in normal coordinates as in Proposition \ref{invariant pair}.  Define $m$ for $M$ and the corresponding  $\hat{m}$ for $\hat{M}$.  If $M$ and $\hat{M}$ are formally equivalent, then $m =  \hat{m}$.
\end{cor}

\begin{proof}
Lemma \ref{m' to m} implies $m \ge \hat{m}$. Reversing the roles of $M$ and $\hat{M}$ yields the other inequality.
\end{proof}

We shall be primarily interested in formal real hypersurfaces which are of infinite type, but nonflat, at $0$.  That is, formal hypersurfaces of $m$-infinite type for some positive integer $m$.  In this case, Corollary \ref{m as invariant} may be strengthened as follows.

\begin{prop} \label{f0 equation found}
If $M$ is of $m$-infinite type at $0$ and $H \in {\F}(M,0;\hat{M},0)$, then $\hat{M}$ is of $m$-infinite type at $0$, $g_0, g_1, \dots, g_{m-1}$ are constant, and
\[
   0 \not \equiv \Theta_{s^m}(z,\chi,0)
   \equiv g_0(0)^{m-1} \, \hat{\Theta}_{{\hat{s}}^m} \big( f_0(z), \bar{f_0}(\chi), 0 \big).
\]
\end{prop}

\begin{proof}
This follows immediately from Lemma \ref{m' to m}, Corollary \ref{m as invariant}, and Lemma \ref{m, equivalently}. 
\end{proof}

We now have the necessary ingredients to prove Proposition \ref{invariant pair}.

\begin{proof}
We have seen that $m = \hat{m}$.  If $m = \hat{m} = 0$, i.e.~the hypersurfaces are of finite type, then it is well known that the triple $(r,L,K)$ is a formal invariant.  (An outline of the proof that $r$ is a formal invariant, for example, may be found in \cite{BaouendiEbenfeltRothschild:RealSubmanifolds}, Chapter I.)  Similarly, observe that $r = \infty$ if and only if $m = \hat{m} = \infty$, which in turn holds if and only if $\hat{r} = \infty$; similarly if $L = \infty$ or $K = \infty$.

Hence, it suffices to assume that all the numbers in question are positive integers. By Proposition \ref{f0 equation found}, we have
\[
   0 \not \equiv \Theta_{s^m}(z,\chi,0)
   \equiv g_0(0)^{m-1} \, \hat{\Theta}_{{\hat{s}}^m} \big( f_0(z), \bar{f_0}(\chi), 0 \big).
\]
A straightforward induction using Faa de Bruno's formula implies that for any multi-indices $\alpha$ and $\beta$,
\begin{align*}
    \Theta_{z^\alpha \chi^\beta s^m}(z,\chi,0) 
    &= g_0(0)^{m-1} 
    \sum_{ \substack{ 
       |\mu| \le |\alpha| \\ 
       |\nu| \le |\beta| } } 
   \hat{\Theta}_{\hat{z}^\mu \hat{\chi}^{\nu} \hat{s}^m} \big( f_0(z), \bar{f_0}(\chi), 0 \big) \\
    &\qquad \qquad \times P^{\alpha \beta}_{\mu \nu} \big( \big( (f^u_0)_{z^\gamma}(z) \big)_{|\gamma| \le |\mu|}, \big( (\bar{f^u_0})_{\chi^\delta}(\chi) \big)_{|\delta| \le |\nu|} \big)
\end{align*}
where each $P^{\alpha \beta}_{\mu \nu}$ is a polynomial in its arguments.

In particular, this implies whenever $|\alpha| + |\beta| < \hat{r}$, we have $\Theta_{z^\alpha \chi^\beta s^m}(0,0,0)= 0$, whence $r \ge \hat{r}$ necessarily.  Reversing the roles of $M$ and $\hat{M}$ implies $r = \hat{r}$.  Similarly, this implies that $\Theta_{\chi^\beta s^m}(z,0,0) \equiv 0$ whenever $|\beta| < \hat{L}$, whence $L \ge \hat{L}$; reversing the roles of the formal hypersurfaces establishes equality.  The proof that $K = \hat{K}$ is similar, and is left to the reader.
\end{proof}


\section{The 1-infinite type case in $\C^2$}

\subsection{Notation and results}

For the remainder of the paper, we shall deal only with formal real hypersurfaces of $\C^2$, and in particular, those hypersurfaces which are of $1$-infinite type at $0$.  Suppose that $M$ is such a formal hypersurface.  We shall write $M$ in normal coordinates $Z=(z,w)$ as in (\ref{eqn:I in Q}).  Since $M$ is of $1$-infinite type, this implies that we can write $Q(z,\chi,\tau) = \tau \, S(z,\chi,\tau)$ for some $S \in \C[[z,\chi,\tau]]$, so that
\begin{equation} \label{eqn:I}
   M =  \bigg\{ \bigg( \frac{w - \bar{w}}{2 i} \bigg) = \Theta \bigg( z, \bar{z}, \frac{w + \bar{w}}{2} \bigg) \bigg\}
     =  \big\{ w = \bar{w} \, S(z,\bar{z},\bar{w}) \big\}.
\end{equation}
For convenience, we shall write
\begin{equation} \label{eqn:p_j defined}
   \theta(z,\chi) = \sum_{j=0}^\infty \frac{\theta_j(z)}{j!} \, \chi^j := \Theta_s(z,\chi,0)
   \not\equiv 0
\end{equation}
Observe that $\theta_j(z) \equiv 0$ if $j < L$ and $\theta_L^{(j)}(0) = 0$ if $j < K$,where $L,K$ are defined by equations (\ref{eqn:L defined}) and (\ref{eqn:K defined}).  It will be useful for later computations to observe that Lemma \ref{m, equivalently} implies
\begin{equation} \label{eqn:S subs 1}
   S(z,\chi,0) = \frac{ 1 + i \, \theta(z,\chi) }{ 1 - i \, \theta(z,\chi) },
\end{equation}
whence repeated differentiation in $\chi$ yields
\begin{equation} \label{eqn:S subs 2}
   S_{\chi^j}(z,0,0) =
   \left\{ \begin{array}{l @{\quad} l}
      1 & j = 0 \\
      0 & 1 \le j \le L - 1 \\
      2i \, \theta_L(z) & j = L \\
      2i \, \theta_{L+1}(z) - 4 \, \theta_1(z)^2 & j = L + 1
   \end{array} \right. .
\end{equation}

We define a new, rather technical, invariant for 1-infinite type hypersurfaces.  Letting $\delta^j_k$ denote the Kronecker delta function (i.e.~$\delta^j_k = 0$ if $j \neq k$, and $\delta^j_j = 1$), we shall define the number $T \in \{0,1\}$ by
\begin{equation} \label{eqn:T defined}
   T := \prod_{q=0}^{K-2} \delta^0_{\theta_{L+1}^{(q)}(0)}.
\end{equation}
That is, $T = 1$ if and only if $\theta_{L+1}(z) = O(|z|^{K-1})$; by means similar to the proofs for the numbers $r$, $L$, and $K$, it can be shown that $T$ is a formal invariant.  Details are left to the reader.

Assume now that $\hat{M}$ is a formal real hypersurface of $\C^2$ which is formally equivalent to $M$, and write it in normal coordinates $\hat{Z} = (\hat{z},\hat{w})$ as
\begin{equation} \label{eqn:I'}
   \hat{M} =  \bigg\{ \frac{\hat{w} - \bar{\hat{w}}}{2 i} = \hat{\Theta} \bigg( \hat{z} , \bar{\hat{z}}, \frac{\hat{w} + \bar{\hat{w}}}{2} \bigg) \bigg\}
           =  \big\{ \hat{w} = \bar{\hat{w}} \, \hat{S}(\hat{z}, \bar{\hat{z}}, \bar{\hat{w}}) \big\},
\end{equation}
Let us write $\hat{\theta}(\hat{z},\hat{\chi}):= \hat{\Theta}_{\hat{s}}(\hat{z},\hat{\chi},0)$ as above.

If $H : (M,0) \to (\hat{M},0)$ is a formal equivalence, then Lemma \ref{form of H} implies that $H(z,w)$ is of the form given by  (\ref{eqn:H expanded}), with $f, g \in \C[[z,w]]$ and $f_z(0,0) \, g(0,0) \neq 0$.  Observe that identity (\ref{eqn:formal map, Q}) can be rewritten (after canceling an extra $\tau$ from both sides) as the identity 
\begin{equation} \label{eqn:formal map, S}
   S(z,\chi,\tau) \, g \big( z, \tau \, S(z,\chi,\tau) \big) 
     \equiv \bar{g}(\chi,\tau) \, \hat{S} \big( f \big( z, \tau \, S(z,\chi,\tau) \big), \bar{f}(z,\chi), \tau \, \bar{g}(\chi,\tau) \big). 
\end{equation}

We shall continue to use the formal Taylor expansions of $f$ and $g$ in $w$ given by equation (\ref{eqn:fn and gn expanded}), and shall write
\begin{equation} \label{eqn:barred derivatives defined}
   f_n(z) := \sum_{k \ge 0} \frac{1}{k!} \, \bar{ a_n^k } \, z^k, \quad 
   g_n(z) := \sum_{k \ge 0} \frac{1}{k!} \, \bar{ b^k_n } \, z^k,
\end{equation}
where the bar denotes complex conjugation.  Note that, in particular,  $a_0^0 = 0$, $a_0^1 \neq 0$, and $b_0^0 = \bar{b_0^0} \neq 0$.

Finally, for each $n \ge 0$, define the formal rational mapping $\Upsilon^n : (\C^2,0) \to (\C^4,0)$ by
\begin{align}
  &\Upsilon^n_1(z,\chi) :=
        K \, \frac{\theta_L(z)}{{\theta_L}'(z)} \bigg(  \frac{ 1 + i\, \theta(z,\chi) }{ 1 - i\, \theta(z,\chi) } \bigg)^n \theta_z(z,\chi)
        - L \, \frac{\bar{\theta_L}(\chi)}{\bar{\theta_L}'(\chi)} \theta_\chi(z,\chi), \label{eqn:UP1}\\
  &\Upsilon^n_2(z,\chi) :=
         (1 + \theta(z,\chi)^2) \bigg[ \bigg(  \frac{ 1 + i\, \theta(z,\chi) }{ 1 - i\, \theta(z,\chi) } \bigg)^n \! \! - 1 \bigg]
        - 2 i \, n \, \frac{\bar{\theta_L}(\chi)}{\bar{\theta_L}'(\chi)} \theta_\chi(z,\chi), \label{eqn:UP2}\\
  &\Upsilon^n_3(z,\chi) :=
        \delta^1_L \, \delta^1_T \bigg\{ \delta^1_K \bigg( \theta_1^{(L)}(0) \frac{\theta_\chi(z,\chi,0)}{\bar{\theta_1}'(\chi)} \bigg) \label{eqn:UP3} \\
  &\qquad  + \bigg( \frac{ \theta_1^{(K)}(0) \theta_2^{(K)}(0) - \theta_1^{(K+1)}(0) \theta_2^{(K-1)}(0) }{ K \, \theta_1^{(K)}(0)^2 } \bigg) \frac{ \bar{\theta_1}(\chi)}{\bar{\theta_1}'(\chi)} \theta_\chi(z,\chi)  \nonumber \\
  &\qquad  - \bigg( \frac{1 + i\,\theta(z,\chi)}{1 - i\, \theta(z,\chi)} \bigg)^n \bigg[ \theta_1(z) \big( 1 + \theta(z,\chi)^2 \big) + \bigg( \frac{\theta_2(z)}{{\theta_1}'(z)} - 2 i \, n \, \frac{\theta_1(z)^2}{{\theta_1}'(z)} \bigg) \theta_z(z,\chi) \bigg]  \nonumber \\
  &\qquad  + \frac{ \theta_2^{(K-1)}(0) }{ \theta_1^{(K)}(0) } \bigg[ \bar{\theta_1}(\chi) \big( 1 + \theta(z,\chi)^2 \big) + \bigg( \frac{\bar{\theta_2}(\chi)}{\bar{\theta_1}'(\chi)} + 2 i \, n \, \frac{\bar{\theta_1}(\chi)^2}{\bar{\theta_1}'(z)} \bigg) \theta_\chi(z,\chi) \bigg] \bigg\}, \nonumber \\
  &\Upsilon^n_4(z,\chi) :=
         \delta^1_K \bigg\{  \frac{ \bar{\theta_1}(\chi) }{{\theta_1}'(0)} (1 + \theta(z,\chi)^2) - \frac{\theta_z(z,\chi)}{{\theta_1}'(z)} \bigg(  \frac{ 1 + i\, \theta(z,\chi) }{ 1 - i\, \theta(z,\chi) } \bigg)^n \label{eqn:UP4}\\
  &\qquad  + \frac{\theta_\chi(z,\chi)}{{\theta_1}'(0)} \bigg[ 2i \, n \, \frac{\bar{\theta_1}(\chi)^2}{ \bar{\theta_1}'(\chi)} + \frac{\bar{\theta_2}(\chi)}{\bar{\theta_1}'(\chi)} - \bigg( \frac{{\theta_1}''(0)}{{\theta_1}'(0)} \bigg) \frac{\bar{\theta_1}(\chi)}{ \bar{\theta_1}'(\chi)} \bigg] \bigg\}, \nonumber
\end{align}
where the $\theta_j$ are defined by equation (\ref{eqn:p_j defined}).  We shall prove in the next chapter that equations (\ref{eqn:UP1}) through (\ref{eqn:UP4}) actually define \emph{formal power series} in $(z,\chi)$, rather than quotients of formal power series.

Observe that the formal mapping $\Upsilon^n$ depends on the choice of normal coordinates $Z=(z,w)$ for the formal hypersurface $M$.

We are now in a position to state the main technical result of the paper, which may be viewed as a sharper version of Theorem \ref{Parametrization}, but with conjugated derivatives.

\begin{thm} \label{linear independence}
Let $(M,0)$ be a formal real hypersurface in $\C^2$ which is of $1$-infinite type, given in normal coordinates $Z = (z,w)$  by equation (\ref{eqn:I}).  Define $\Upsilon^n(z,\chi)$ by equations (\ref{eqn:UP1}) through (\ref{eqn:UP4}).  For each $n \in \N$, define the complex vector space
\begin{equation} \label{eqn:V defined}
   \V^n(M) := \Span_\C \setof{ \upsilon^n_{s,t}:= \Upsilon^n_{z^s \chi^t}(0,0) }{ s,t \in \N } \subset \C^4.
\end{equation}
Then the dimension of the vector space $\V^n(M)$ is a formal invariant for each $n$, and the invariant set of integers
\begin{equation} \label{eqn:D defined}
   \D(M) := \setof{ n \in \N }{ \dim_\C \V^n(M) < 2 + \delta^1_K + \delta^1_L \, \delta^1_T } 
\end{equation}
is always finite.

Furthermore, given any formal real hypersurface $(\hat{M},0)$ in $\C^2$ formally equivalent to $(M,0)$, any normal coordinates $\hat{Z} = (\hat{z},\hat{w})$ for $\hat{M}$, and any $n \in \N$, there exists a formal power series $\A_n(z;\Delta, \Lambda) \in \C[\Delta,\Lambda][[z]]^2$, with $(z,\Delta, \Lambda) \in \C \times \C \times \C^{4 |\D(M)|}$, such that 
\[
   \big( f_n(z), g_n(z) \big) \equiv \A_n \bigg( z ; \frac{1}{a_0^1 b_0^0}, \big( a_j^0, b_j^0, a_j^1, b_j^1 \big)_{j \in \D(M)} \bigg).
\]
for any  $H \in \F(M,0;\hat{M},0)$.

Moreover, if $M$ and $\hat{M}$ are convergent, then there exists an $\epsilon > 0$ such that 
\[
   \bigg\{ z \mapsto \A_n \bigg( z ; \frac{1}{a_0^1 b_0^0} , \big( a_j^0, b_j^0, a_j^1, b_j^1 \big)_{j \in \D(M)} \bigg) \bigg\} \in \calO_\epsilon(z)^2
\]
for every $H \in \F(M,0;\hat{M},0)$ and every $n \in \N$.
\end{thm}

\subsection{Examples}

In this section, we use Theorem \ref{linear independence} and Proposition \ref{f0 equation found} to calculate the formal transformation groups of various $1$-infinite type hypersurfaces.

\begin{example}
Consider the family of $1$-infinite type hypersurfaces
\[
   M_c^j := \setof{ (z,w) }{ \Im \, w = c (\Re \, w) |z|^{2j} }, \quad c \in \R \setminus \{0\}, \quad j \ge 1.
\]
Observe that $L = K = j$, $T = 1$, and $\theta(z,\chi) = c \, z \chi$.  If $n > 0$, it can be shown that $\{ \upsilon^n_{2j,2j}, \upsilon^n_{3j,3j} \}$ is a basis for $\V^n(M_c^j)$ if $j \ge 2$, and that adding the vectors $\{ \upsilon^n_{2,3}, \upsilon^n_{3,2} \}$ extends this to a basis for $\V^n(M_c^1)$.  Hence, in any case, we have $\D(M_c^j) = \{ 0 \}$, so any formal equivalence with source $M_c^j$ is determined by $(a_0^1,b_0^0)$.

Applying Proposition \ref{f0 equation found} with $M = \hat{M} = M^j_c$ implies $f_0(z) = \varepsilon \, z$ for some $\varepsilon \in \C$ with $|\varepsilon| = 1$.  It thus follows that
\[
   {\Aut}(M_c^j,0) = \setof{ (z,w) \mapsto \big( \varepsilon \, z, r \, w \big) }{ 
   \varepsilon \in \C, \, |\varepsilon| = 1, \, r \in \R \setminus \{ 0 \} }.
\]
In particular, every formal automorphism converges.

Observe that for $j \neq k$, the hypersurfaces $M_c^j$ and $M_b^k$ are \emph{not} formally equivalent (Theorem \ref{4-tuple}).  On the other hand, $M^j_c$ and $M^j_b$ are formally equivalent if and only if $c / b > 0$.  In this case, applying Proposition \ref{f0 equation found} implies that $f_0(z) = \alpha \, z$ for some $\alpha \in \C$ of modulus $(c/b)^{1/{2j}}$.  It thus follows that
\[
   {\F}(M_c^j,0;M_b^j,0) =  \setof{ (z,w) \mapsto \bigg( \frac{c}{b} \bigg)^{ \frac{1}{2j} } \big( \varepsilon z, r \, w \big) }{ 
   \varepsilon \in \C, \, |\varepsilon| = 1, \, r \in \R \setminus \{ 0 \} }.
\]
Hence, the hypersurfaces $M_c^j$ are formally equivalent if and only if they are biholomorphically equivalent if and only if $b$ and $c$ have the same sign.
\end{example}

\begin{example}
Consider the family of $1$-infinite type hypersurfaces
\[
   N_b^j := \setof{ (z,w) }{ \Im \, w = 2 (\Re \, w) (\Re(b \, z \bar{z}^j) }, \quad b \in \C \setminus \{ 0\}, \, j \ge 2.
\]
Note $L = 1$, $K = j$, and $\theta(z,\chi) = b \, z \chi^j + \bar{b} \, z^j \chi$. If $n > 0$, it can be shown that $\{ \upsilon^n_{2,2}, \upsilon^n_{3,2}, \upsilon^n_{3,3} \big\}$ forms a basis for $\V^n(N_b^j)$, so we again conclude that $\D(N_b^j) = \{ 0 \}$.  Hence, every formal equivalence  $H$ with source $N_b^j$ is determined by the values $a_0^1$ and $b_0^0$.

Now, Proposition \ref{f0 equation found} applied to the $M = \hat{M} = N_b^j$ case implies that $a_0^1$ is a $(j-1)$-th root of unity, and that $f_0(z) = z / a_0^1$. We conclude
\[
   {\Aut}(N_b^j,0) = \setof{ (z,w) \mapsto \big( \varepsilon \, z, r \, w \big) }{ 
      \varepsilon \in \C,  \, \varepsilon^{j-1} = 1, \, r \in \R \setminus \{0\} }.
\]
Note that every formal automorphism converges.
\end{example}

\begin{example}
Consider the hypersurface
\[
    B_0 := \setof{ (z,w) }{ \Im \, w = (\Re \, w) \frac{1 - \sqrt{1 - 4 z^2 \chi^2}}{2 z \chi} }.
\]
It is easy to check that $L = K = 1$ in this case, and that $\D(B_0) = \{ 0, 1, 2 \}$.   (In fact, we have $\Upsilon^1_4 \equiv 0$, and $2i \, \Upsilon^2_1 \equiv \Upsilon^2_2$.)  The author has calculated the entire stability group of the hypersurface $B_0$, which is an example of a real-analytic hypersurface whose stability group at the origin is determined by $3$-jets \emph{but not by $2$-jets}; see \cite{Kowalski:AHypersurfaceWhoseStabilityGroup}.

In general, for any integer $n > 0$, there exists a (unique) real-valued power series $\rho_n(t)$ with $\rho_n(0)=0$ and ${\rho_n}'(0) = 1$, such that for the $1$-infinite type 2 hypersurface
\[
   B_{n}:= \setof{ (z,w) }{ \Im \, w = (\Re \, w) \rho_n(z \chi) },
\]
we have $\Upsilon^{n}_3 \equiv 0$, and so $n \in \D(B_n)$ necessarily.  That is, while $\D(M)$ always contains only finitely many integers, the integers themselves can be arbitrarily large.  Further examples may be found in \cite{Kowalski:THESIS}, Chapter 7.
\end{example}


\section{Proofs of the main results}
\label{Proofs of the Main Results}

\subsection{Proof of Theorem \ref{linear independence}}

A basic outline of the proof can be divided into four steps.
\begin{enumerate}
  \item Given a fixed set of normal coordinates $Z = (z,w)$, we prove that for each $n \in \N$ the power series $f_n(z)$ and $g_n(z)$ are rationally parametrized by the values $(a_\ell^j, b_\ell^j)$ for $\ell = 0,1$ and $0 \le j \le n$.  
  \item We prove that under these conditions, if $n \not \in \D(M)$, then the $4$-tuple of complex numbers $(a_n^0,a_n^1,b_n^0,b_n^1)$ is itself a polynomial in $1/(a_0^1 \, b_0^0)$ and $(a_\ell^j, b_\ell^j)$ for $\ell = 0,1$ and $0 \le j \le n - 1$.  
  \item We prove that $\D(M)$, defined by these normal coordinates, is always finite.  
  \item We show that the dimension of $\V^n(M)$ (and hence the set $\D(M)$) is independent of the normal coordinates used to define it.  
\end{enumerate}

To fix notation throughout the proof, we shall assume that $M$ is always given in normal coordinates $Z = (z,w)$ by (\ref{eqn:I}).  We shall also set $\D = \D(M)$ and $\V^n = \V^n(M)$.  Similarly, $\hat{M}$, whenever a target formal hypersurface is needed, will always be given in normal coordinates $\hat{Z}=(\hat{z},\hat{w})$ by (\ref{eqn:I'}).  If $H : (M,0) \to (\hat{M},0)$ is a formal equivalence, we shall set
\begin{align*}
   \Delta(H) &:= \frac{ 1}{ a_0^1 b_0^0} \in \C \setminus \{ 0 \}, \\
   \lambda^n_2(H) &:= \big( a_n^1, b_n^0 \big) \in \C^2, \\
   \lambda^n_3(H) &:= \big( a_n^1, b_n^0, a_n^0 \big) \in \C^3,\\
   \lambda^n_4(H) &:= \big( a_n^1, b_n^0, a_n^0, b_n^1 \big) \in \C^4, \\
   \Lambda^n_j(H) &:= \big( \lambda^0_j(H), \lambda^1_j(H), \dots, \lambda^n_j(H) \big) \in \C^{j(n+1)}.
\end{align*}

We shall also use the following conventions for naming various types of polynomials and power series.
\begin{itemize}
  \item $\calQ^d(X;\Lambda) \in \C[X,\Lambda] \equiv \C[\Lambda] [X]$ denotes a polynomial in $X$ of degree $d$ whose coefficients are polynomial in $\Lambda$.
  \item $\calP(\Lambda;X) \in \C[[X,\Lambda]] \equiv \C[[X]] [\Lambda]$ denotes a polynomial in $\Lambda$ whose coefficients are power series in $X$.
  \item $\calR(X;\Lambda) \in \C[[X,\Lambda]] \equiv \C[\Lambda] [[X]]$ denotes a power series in $X$ whose coefficients are polynomial in $\Lambda$.
\end{itemize}

Let us assume the normal coordinates $Z$ and $\hat{Z}$ for $M$ and $\hat{M}$ are fixed.  We now tackle the first step, the parametrizing of $f_n$ and $g_n$.  We begin with a lemma.

\begin{lem} \label{reflection identities 0}
Let $(M,0)$ and $(\hat{M},0)$ be formally equivalent formal $1$-infinite type hypersurfaces as above.  Then there exist unique formal power series $U,V \in \C[[X,Y]]$, vanishing at $0$, such that
\[
   f_0(z) = U \bigg( z, \frac{z}{a_0^1} \bigg), \quad \bar{f_0}(\chi) = V \bigg( \chi, a_0^1 \, \chi \bigg)
\]
for any $H \in {\F}(M,0;\hat{M},0)$.  Moreover, if both $M$ and $\hat{M}$ are convergent hypersurfaces, then $U,V \in \C\{X,Y\}$.
\end{lem}

\begin{proof}
Proposition \ref{f0 equation found} implies that
\begin{equation} \label{eqn:f0 in phi}
   \theta( z, \chi) \equiv \hat{\theta} \big( f_0(z), \bar{f_0}(\chi) \big).
\end{equation}
Differentiating this $L$ times in $\chi$ using Faa de Bruno's formula and setting $\chi = 0$ yields the identity
\begin{equation} \label{eqn:f0 fuller}
   \theta_L(z) \equiv (a_0^1)^L \, \hat{\theta}_L \big( f_0(z) \big).
\end{equation}
Differentiating this $K$ times in $z$ and setting $z = 0$ yields
\begin{equation} \label{eqn:plk elimination}
   \theta_L^{(K)}(0) =  \big( \, \bar{a_0^1} \, \big)^K \, (a_0^1)^L \,  \hat{\theta}_L^{(K)}(0).
\end{equation}
In particular, we find that for any formal equivalence $H \in {\F}(M,0;\hat{M},0)$,
\begin{equation} \label{eqn:a01 constant norm}
   \big| {f_0}'(0) \big| = | a_0^1 | = \left| \frac{ \theta_L^{(K)}(0) }{ \hat{\theta}_L^{(K)}(0) }\right|^{ \frac{1}{L+K} } =: \mu \in \R \setminus \{ 0 \}.
\end{equation}

Now, observe we can write
\[
   \theta_L(z) = \frac{1}{K!} \theta_L^{(K)}(0) \, z^K \, t(z),
\]
for some $t \in \C[[z]]$ with $t(0) = 1$.  Thus, there exists a unique power series $u(z)$ with $u(0) = 1$ such that $u(z)^K = t(z)$.  Similarly, let us write
\[
   \hat{\theta}_L(\hat{z}) = \frac{1}{K!} \hat{\theta}_L^{(K)}(0) \, \hat{z}^K \, \hat{u}(\hat{z})^K,
\]
with $\hat{u}(0)=1$. Define the formal power series
\[
   \iota(\hat{z},X,Y) := \hat{z} \, \hat{u}(\hat{z}) - \mu^2 \,  Y \, u(X) .
\]
Observe that $\iota(0,0,0) = 0$ and $\iota_{\hat{z}}(0,0,0) = 1$, whence the formal Implicit Function Theorem implies the existence of a unique power series $U(X,Y)$, vanishing at $(0,0)$, such that $\iota \big( U(X,Y),X,Y \big) \equiv 0$.

Now, suppose that $H \in {\F}(M,0;\hat{M},0)$.  Then identity (\ref{eqn:f0 fuller}) may be written as
\[
   \frac{1}{K!} \theta_L^{(K)}(0) \, \big( z \, u(z) \big)^K \equiv 
   (a_0^1)^L \frac{1}{K!} \hat{\theta}_L^{(K)}(0) \, \bigg( f_0(z) \, \hat{u}\big( f_0(z) \big) \bigg)^K.
\]
Replacing ${\theta_L}^{(K)}(0)$ by equation (\ref{eqn:plk elimination}) and canceling common terms yields the identity
\[
   \bigg[ \, \bar{a_0^1} \, z \, u(z) \bigg]^K \equiv \bigg[ f_0(z) \, \hat{u} \big( f_0(z) \big) \bigg]^K.
\]
Formally extracting $K$-th roots on both sides, we conclude that the two power series in the brackets differ only by some multiple $\varepsilon \in \C$ with $\varepsilon^K = 1$.  However, since
\[
   \dd{}{z} \bigg\{\, \bar{a_0^1} \, z \, u(z) \bigg\} \bigg|_{z=0}
   = \bar{a_0^1} = {f_0}'(0) 
   = \dd{}{z} \bigg\{ f_0(z) \, \hat{u} \big( f_0(z) \big) \bigg\} \bigg|_{z=0},
\]
we conclude that $\varepsilon = 1$ necessarily.  Moreover, since $a_0^1 \, \bar{a_0^1} = \mu^2$, we have
\[
   \mu^2 \bigg( \frac{z}{a_0^1} \bigg)  u(z) \equiv  f_0(z) \, \hat{u} \big( f_0(z) \big).
\]    
Hence, $\iota \big( f_0(z), z, z/a_0^1) \equiv 0$, so by the uniqueness of $U$, we conclude $f_0(z) = U \big( z, z / a_0^1 \big)$.  Conjugating this result yields $\bar{f_0}(\chi) = V(\chi, a_0^1 \, \chi)$, where $V$ is defined by $V(X, Y) := \bar{U}( X, Y/\mu^2 )$.

Finally, observe that if $M$ and $\hat{M}$ are convergent, then the power series $\theta$ (hence $u$) and $\hat{\theta}$ (hence $\hat{u}$) are convergent, so the holomorphic Implicit Function Theorem implies that $U$ and $V$ are necessarily convergent near $(0,0) \in \C^2$.
\end{proof}

We can now extend this lemma to show that $f_n$ and $g_n$ are similarly parametrized for any $n \ge 0$.

\begin{prop} \label{reflection identities}
Let $(M,0)$, $(\hat{M},0)$ be formally equivalent formal $1$-infinite type hypersurfaces as above.  Then for every $n \in \N$, there exists a power series $\B_n(z;\Delta,\Lambda) \in \C[\Delta,\Lambda][[z]]^2$ such that the following holds for any $H \in {\F}(M,0;\hat{M},0)$:
\begin{equation}  \label{eqn:finite det, sort of}
   \big( f_n(z), g_n(z) \big) = \B_n \bigg( z; \Delta(H), \Lambda^n_{2+\delta^1_K+\delta^1_T}(H) \bigg).
\end{equation}
In addition, if $n \ge 1$, then in fact
\begin{align}
\label{eqn:fn equation formed}
   &\frac{f_n(z)}{{f_0}'(z)} =
         T_n^1 \big( z; \Delta(H), \Lambda^{n-1}_{2+\delta^1_K+\delta^1_T}(H) \big) -  \frac{L}{a_0^1} \bigg[ \frac{\theta_L(z)}{{\theta_L}'(z)} \bigg] a_n^1
         + \frac{n}{b_0^0} \bigg[ \frac{\theta_L(z)}{{\theta_L}'(z)} \bigg]  b_n^0\\
     &\quad  
         + \frac{i \, \delta^1_K}{2 \, b_0^0} \bigg[ \frac{1}{{\theta_1}'(z)} \bigg] b_n^1
         + \frac{\delta^1_T}{a_0^1} \bigg[ 2i \, n \,  \, \frac{\theta_1(z)^2}{{\theta_L}'(z)} 
         -   \frac{\theta_{L+1}(z)}{{\theta_L}'(z)} + \frac{L \, a_0^2}{ a_0^1} \frac{\theta_L(z)}{{\theta_L}'(z)} \bigg] a_n^0 \nonumber  \\
\label{eqn:gn equation formed}
   &g_n(z) = T_n^2 \big( z; \Delta(H), \Lambda^{n-1}_{2+\delta^1_K+\delta^1_T}(H) \big) +  b_n^0 +  \frac{2 i \, b_0^0 \, \delta^1_T}{a_0^1} \big[ \theta_1(z) \big] a_n^0
\end{align}
with $T(z;\Delta,\Lambda^{n-1}_{2+\delta^1_K+\delta^1_T}) \in \C^2[\Delta,\Lambda^{n-1}_{2+\delta^1_K+\delta^1_T}][[z]]$.

Moreover, if $M$ and $\hat{M}$ are convergent, then there exists an $\epsilon > 0$ such that
\[
    \bigg\{ z \mapsto \B_n \bigg( z; \Delta(H), \Lambda^n_{2+\delta^1_K+\delta^1_T}(H) \bigg) \bigg\}
    \in \calO_\epsilon(z)^2
\]
for every $n \in \N$ and every $H \in {\F}(M,0;\hat{M},0)$.
\end{prop}

\begin{proof}
For convenience, we shall set  $\gamma = 2 + \delta^1_K + \delta^1_T$.  We proceed by induction.  The $n=0$ case follows immediately from Lemma \ref{reflection identities 0} and the fact that $g_0(z) \equiv b_0^0$  (Proposition \ref{f0 equation found}), so let us assume that the proposition is true up to some $n-1 \ge 0$.  To prove (\ref{eqn:finite det, sort of}), it suffices to prove that equations (\ref{eqn:fn equation formed}) and (\ref{eqn:gn equation formed}) hold.

Suppose that $H : (M,0) \to (\hat{M},0)$ is a formal equivalence.\footnote{We remark that the construction given in this section can be carried out if \emph{no} formal equivalence exists between $M$ and $\hat{M}$.}  Differentiating identity (\ref{eqn:formal map, S}) $n$ times in $\tau$ using {Faa de Bruno's formula} and setting $\tau = 0$ (or, equivalently, substituting $Q(z,\chi,\tau) = \tau \, S(z,\chi,\tau)$ and $v = n+1$ into identity (\ref{eqn:formal map, Q tau})) yields
\begin{multline} \label{eqn:formal map in tau}
   - S(z,\chi,0)^{n+1} g_n(z) + b_0^0 \, \hat{S}_{\hat{z}} \big( f_0(z), \bar{f_0}(\chi), 0 \big) \, S(z,\chi,0)^n \, f_n(z) \\
    + b_0^0 \, \hat{S}_{\hat{\chi}} \big( f_0(z), \bar{f_0}(\chi), 0 \big)  \, \bar{f_n}(\chi) + \hat{S} \big( f_0(z), \bar{f_0}(\chi), 0 \big) \bar{g_n}(\chi) \\
    \equiv \calP_n \bigg( b_0^0,  \big( f_j(z), g_j(z), \bar{f_j}(\chi), \bar{g_j}(\chi) \big)_{j=1}^{n-1}; z, \chi, f_0(z), \bar{f_0}(\chi) \bigg),
\end{multline}
where $\calP_n(\Lambda;X)$, with $(\Lambda,X) \in \C^{4n-3} \times \C^4$, depends only on $M$ and $\hat{M}$ and \emph{not} the map $H$.\footnote{Indeed, an explicit formula for $\calP_n$ is given following the proof of Proposition \ref{reflection identities}.}  Note that Lemma \ref{m' to m} implies $\hat{S} \big( f_0(z), \bar{f_0}(\chi),0 \big) = S(z,\chi,0)$, whence
\[
   \hat{S}_{\hat{z}} \big( f_0(z), \bar{f_0}(\chi),0 \big) = \frac{ S_z(z,\chi,0) }{ {f_0}'(z) }, \quad
   \hat{S}_{\hat{\chi}} \big( f_0(z), \bar{f_0}(\chi),0 \big) = \frac{ S_\chi(z,\chi,0) }{ \bar{f_0}'(\chi) }
\]

Observe that if equation (\ref{eqn:finite det, sort of}) holds for some $n \in \N$, then
\begin{equation} \label{eqn:barring lambda}
   \bar{\lambda^n_4(H)} = \big( (\B_n)^1_z, (\B_n)^2, 
       (\B_n)^1, (\B_n)^2_z \big)(0;\Delta(H),\Lambda^n_\gamma(H)) 
    =: \beta_n(\Delta(H),\Lambda^n_\gamma(H)).
\end{equation}
Applying the inductive hypothesis to this and substituting this into equation (\ref{eqn:formal map in tau}) yields
\begin{equation} \label{eqn:barring the fn and gn's}
   \big( \, \bar{f_j}(\chi), \bar{g_j}(\chi) \big) 
    = \bar{\B_j} \bigg( \chi; \bigg( \frac{ a_0^1 }{\mu} \bigg)^2 \Delta(H), \big( \beta_\ell(\Delta(H),\Lambda^\ell_\gamma(H)) \big)_{\ell=0}^j \bigg)
\end{equation}
for $j < n$, where $\mu$ is defined in equation (\ref{eqn:a01 constant norm}).  Substituting these values into  (\ref{eqn:formal map in tau}) yields
\begin{multline}\label{eqn:n pair expression}
   - S(z,\chi,0)^{n+1} g_n(z)  + S(z,\chi,0) \bar{g_n}(\chi)  + b_0^0 \,  S_z(z,\chi,0)  S(z,\chi,0)^n \, \frac{ f_n(z) }{ {f_0}'(z) }  \\
    + b_0^0 \,  S_\chi(z,\chi,0) \, \frac{ \bar{f_n}(\chi)  }{ \bar{f_0}'(\chi) } 
    \equiv \calR_n(z,\chi;\Delta(H), \Lambda^{n-1}_\gamma(H) ),
\end{multline}
with $\calR_n(X;\Lambda)$ independent of the mapping $H$ for each $n \ge 0$.

On one hand, substituting $\chi = 0$ and the identities from equations (\ref{eqn:S subs 1}) and (\ref{eqn:S subs 2}) into (\ref{eqn:n pair expression}) yields
\begin{equation} \label{eqn:gn(z) defined}
   g_n(z) = \calR_n(z,0;\Delta(H),\Lambda^{n-1}_\gamma(H)) +  b_n^0 +  \frac{2 i \, b_0^0}{a_0^1} \big[ \theta_1(z) \big] a_n^0.
\end{equation}

On the other hand, differentiating identity (\ref{eqn:n pair expression}) $L$ times in $\chi$, setting $\chi = 0$, and using the identities from equations (\ref{eqn:S subs 1}) and (\ref{eqn:S subs 2}) yields (after rearranging terms) the identity
\begin{multline*}
  {\theta_L}'(z) \, \frac{ f_n(z) }{ {f_0}'(z) }
   \equiv - \frac{i}{2 \, b_0^0} (\calR_n)_{\chi^j}(z,0;\Delta(H),\Lambda^{n-1}_\gamma(H)) + \frac{ (n+1) }{ b_0^0 } \theta_L(z) \, g_n(z)  + \frac{i}{2 \, b_0^0} b_n^L\\
    - \frac{1}{b_0^0} \big[ \theta_L(z) \big] b_n^0 
          - \frac{ L }{a_0^1} \big[ \theta_L(z) \big] a_n^1
          - \frac{1}{a_0^1} \bigg[  \theta_{L+1}(z) + 2 i \, \theta_1(z)^2 - \frac{  L \, a_0^2}{a_0^1} \, \theta_L(z) \bigg] a_n^0.
\end{multline*}
Using the formula for $g_n(z)$ from equation (\ref{eqn:gn(z) defined}) and observing that $(\theta_1)^2 = \theta_1 \, \theta_L$ for every $L \ge 1$, we can rewrite this identity as
\begin{multline} \label{eqn:fn(z) defined}
  {\theta_L}'(z) \, \frac{ f_n(z) }{ {f_0}'(z) }
   \equiv - \frac{i}{2 \, b_0^0} (\calR_n)_{\chi^j}(z,0;\Delta(H),\Lambda^{n-1}_\gamma(H)) 
          - \frac{n}{b_0^0} \big[ \theta_L(z) \big] b_n^0 + \frac{i}{2 \, b_0^0} b_n^L  \\
    - \frac{ L }{a_0^1} \big[ \theta_L(z) \big] a_n^1 
          + \frac{1}{a_0^1} \bigg[ - \theta_{L+1}(z) + 2 i \, n \, \theta_1(z)^2 + \frac{  L \, a_0^2}{a_0^1} \, \theta_L(z) \bigg] a_n^0
\end{multline}
We complete the proof by examining cases.

\textbf{Case 1:} $K = 1$.  In this case $L = T =1$ necessarily, so $\gamma = 4$ and ${\theta_L}'(z) = {\theta_1}'(z)$ is a multiplicative unit.  Dividing it on both sides of (\ref{eqn:fn(z) defined}) yields (\ref{eqn:fn equation formed}); equation (\ref{eqn:gn equation formed}) follows from (\ref{eqn:gn(z) defined}).

\textbf{Case 2:} $K > 0$.  In this case, setting $z = 0$ in (\ref{eqn:fn(z) defined}) yields
\[
  0 = - \frac{i}{2 \, b_0^0} (\calR_n)_{\chi^j}(z,0;\Delta(H), \Lambda^{n-1}_\gamma(H)) + \frac{i}{2 \, b_0^0} b_n^L,
\]
whence we may replace $b_n^L$ in identity (\ref{eqn:fn(z) defined}) by $(\calR_n)_{\chi^j}(z,0;\Delta(H), \Lambda^{n-1}_\gamma(H))$.  Thus, after rearranging the terms again, we may rewrite (\ref{eqn:fn(z) defined}) as
\begin{multline} \label{eqn:fn(z) defined 2}
  {\theta_L}'(z) \, \frac{ f_n(z) }{ {f_0}'(z) }
   \equiv \sum_{j=0}^{K-2} \left[ \frac{r^n_j(\Delta(H),\Lambda^{n-1}_\gamma(H))}{j!} z^j + \calR^1_n(z;\Delta(H),\Lambda^{n-1}_0(H)) \right] \\
   - \frac{n}{b_0^0} \big[ \theta_L(z) \big] b_n^0 
   - \frac{ L }{a_0^1} \big[ \theta_L(z) \big] a_n^1
   + \frac{1}{a_0^1} \bigg[ - \theta_{L+1}(z) + 2 i \, n \, \theta_1(z)^2 + \frac{  L \, a_0^2}{a_0^1} \, \theta_L(z) \bigg] a_n^0
\end{multline}
with the $r^n_j$ polynomials and $\calR^1_n(z;\Delta,\Lambda)$ of order at least $K-1$ in $z$.

\textbf{Subcase A:} $T = 1$.  Note that $\gamma = 3$.  Since ${\theta_{L+1}}^{(j)}(0) = 0$ for $j < K-1$, differentiating (\ref{eqn:fn(z) defined 2}) in $z$ (up to $K-2$ times) yields the relations
\[
   r^n_j(\Delta(H), \Lambda^{n-1}_3(H)) = 0, \quad 0 \le j \le K - 2.
\]
Observe that this does not imply that the polynomials $r^n_j(\Delta,\Lambda)$ are themselves identically zero; merely that they vanish whenever 
\[
   (\Delta, \Lambda) = \big( \Delta(H),\Lambda^{n-1}_3(H) \big)
\]
for some formal equivalence $H \in {\F}(M,0;\hat{M},0)$.

Consequently, we may remove the first $K-1$ summands of the right-hand expression in identity (\ref{eqn:fn(z) defined 2}).  Observe that all the remaining summands are of order at least $K-1$ in $z$, and hence can be divided by ${\theta_L}'(z)$ to form another power series.  This division yields (\ref{eqn:fn equation formed}); (\ref{eqn:gn equation formed}) follows from (\ref{eqn:gn(z) defined}).

\textbf{Subcase B:} $T = 0$.  Note that $\gamma = 2$.  We know there exists some $j_0 \in \{ 1, 2, \dots, K-2 \}$ such that ${\theta_{L+1}}^{(j_0)}(0) \neq 0$.  Differentiating the identity (\ref{eqn:fn(z) defined 2}) $j_0$ times in $z$ and setting $z=0$, we obtain
\[
   0 = r^n_{j_0}(\Delta(H),\Lambda^{n-1}_2(H))  - \frac{{\theta_{L+1}}^{(j_0)}(0)}{a_0^1} \, a_n^0,
\]
whence we may replace $a_n^0$ in (\ref{eqn:gn(z) defined}) and (\ref{eqn:fn(z) defined 2}) by $ \frac{ a_0^1 \, r^n_{j_0}(\Delta(H),\Lambda^{n-1}_2(H)) }{ {\theta_{L+1}}^{(j_0)}(0) } $ to obtain
\begin{align*}
  {\theta_L}'(z) \, \frac{ f_n(z) }{ {f_0}'(z) }
   &\equiv \sum_{j=0}^{K-2} \left[ \frac{\tilde{r}^n_j(\Delta(H),\Lambda^{n-1}_2(H))}{j!} z^j + \calR^2_n(z;\Delta(H),\Lambda^{n-1}_2(H)) \right] \\ 
   &\qquad - \frac{n}{b_0^0} \big[ \theta_L(z) \big] b_n^0  
                  - \frac{ L }{a_0^1} \big[ \theta_L(z) \big] a_n^1,\\
  g_n(z) &= \calR^3_n(z,0;\Delta(H),\Lambda^{n-1}_2(H)) +  b_n^0
\end{align*}
Thus, (\ref{eqn:gn equation formed}) holds; arguing as in the proof of Subcase A now yields (\ref{eqn:fn equation formed}).

The only thing missing from the proof is the convergence statement.  Assume now that $M$ and $\hat{M}$ define real-analytic hypersurfaces in $\C^2$ through $0$.  Hence, there exists a $\delta > 0$ such that
\[
   S(z,\chi,\tau) \in \calO_\delta(z,\chi,\tau), \quad
   \hat{S} \big( \hat{z}, \hat{\chi}, \hat{\tau} \big) \in \calO_\delta \big( \hat{z}, \hat{\chi}, \hat{\tau} \big).
\]
Without loss of generality, we shall assume that $\delta$ is chosen small enough such that $\theta_L(z) \neq 0$ for $0 < |z| < \delta$, since the zeros of a nonconstant holomorphic function of one variable are isolated.

Similarly, since $U(X,Y) \in \C\{X,Y\}$ vanishes at 0 by Lemma \ref{reflection identities 0}, there exists an $\eta > 0$ such that $U(X,Y) \in \calO_\eta(X,Y)$ and satisfies 
\[
   \big| U(X,Y) \big| < \delta \quad \mbox{whenever} \quad
   |X|, |Y| < \eta.
\]

Choose $\epsilon < \min \{ \delta, \eta, \mu \, \eta \}$,  where $\mu$ is defined by equation (\ref{eqn:a01 constant norm}).  We claim this is the desired $\epsilon > 0$; the proof is by induction.  The case $n = 0$ follows from Lemma \ref{reflection identities 0}.  Assuming this choice of $\epsilon$ holds up to some $n - 1$, then observe that the mapping
\begin{multline*}
     (z,\chi) \mapsto \calR_n \big( z,\chi;\Delta(H),\Lambda^{n-1}_\gamma(H) \big)\\
    \equiv  \calP_n \bigg( b_0^0,  \big( f_j(z), g_j(z), \bar{f_j}(\chi), \bar{g_j}(\chi) \big)_{j=1}^{n-1}; z, \chi, f_0(z), \bar{f_0}(\chi) \bigg)
\end{multline*}
converges if $|z|, |\chi| < \delta$ for any $H \in {\F}(M,0;\hat{M},0)$. Fix such an $H$.  By equation (\ref{eqn:gn(z) defined}), we conclude $g_n(z)$ converges on the ball $B^1(0,\epsilon) = \setof{ z \in \C }{ |z| < \epsilon }$.  On the other hand, we have shown that
\[
   {\theta_L}'(z) \frac{ f_n(z) }{ {f_0}'(z) } = z^{K-1} q \big( z; \Delta(H), \Lambda^{n-1}_\gamma(H) \big)
\]
with $q(\cdot;\Delta(H),\Lambda^{n-1}_\gamma(H))$ convergent on $B^1(0,\epsilon)$.  Since ${\theta_L}'(z)$ converges for $|z| < \epsilon$ and in the $\epsilon$-ball vanishes only at $z=0$ (of order $K-1$), we conclude that $f_n(z)$ converges on $B^1(0,\epsilon)$ as well, which completes the proof.
\end{proof}

It is of interest to note that as a consequence of  Proposition \ref{reflection identities}, we see that if $M$ and $\hat{M}$ are real-analytic hypersurfaces in $\C^2$ and $H$ is a formal equivalence between them, then the formal mappings $z \mapsto H_{w^n}(z,0)$ are convergent for every $n \in \N$; moreover, they converge on some common $\epsilon$-neighborhood of $0 \in \C$, with $\epsilon$ independent of $n$ and $H$.

Because it is useful in doing calculations, we now give the explicit formula for $\calP_n$.  Using Faa de Bruno's formula, we have
\begin{multline*}
   \calP_n \bigg( \big( f_j,g_j,\bar{f_j},\bar{g_j} \big)_{j=0}^{n-1}; z, \chi, \hat{z}, \hat{\chi} \bigg)\\
   = p_n \bigg( \big( f_j,g_j,\bar{f_j},\bar{g_j} \big)_{0 \le j \le n-1}, 
                \big( S_{\tau^j}(z, \chi,0) \big)_{0 \le j \le n},  
                \big( \hat{S}_{\hat{z}^j \hat{\chi}^k \hat{\tau}^\ell}(\hat{z}, \hat{\chi},0 \big)_{0 \le j + k + \ell \le n} \bigg)
\end{multline*}
where $p_n$ is the universal polynomial
\begin{multline*}
   p_n \bigg( \big( f_j,g_j,\bar{f_j},\bar{g_j} \big)_{0 \le j \le n-1}, 
                \big( S_{j} \big)_{0 \le j \le n},  
                \big( \hat{S}_{(j,k,\ell)} \big)_{0 \le j + k + \ell \le n} \bigg) \\
   \equiv
    \sum_{ \substack{ \alpha \in \N^n \\ k + [\alpha] = n \\ |\alpha| < n } }
     \frac{n! \, g_{|\alpha|} \, S_{k} }{k! \, \alpha!} 
      \prod_{p=1}^n \bigg( \frac{ S_{p-1}}{(p-1)!} \bigg)^{\alpha_p}
    - \sum_{ \substack{ \alpha, \beta, \gamma \in \N^n \\
                     k+[\alpha]+[\beta]+[\gamma]= n \\
                     [\alpha], [\beta], k < n } }
     \frac{n! \, \bar{g_k} \, \hat{S}_{ (|\alpha|, |\beta|, |\gamma|) }}{k! \, \alpha! \, \beta! \, \gamma!} \\
    \times \prod_{p=1}^n \left( \sum_{ \substack{ \xi \in \N^p \\ [\xi]=p } }
           \frac{ f_{|\xi|} }{ \xi! }
             \prod_{q=1}^n \bigg( \frac{ S_{q-1}}{(q-1)!} \bigg)^{\xi_q} \right)^{\alpha_p}
          \bigg( \frac{ \bar{f_p} }{p!} \bigg)^{\beta_p}
          \bigg( \frac{ \bar{g_{p-1}} }{(p-1)!} \bigg)^{\gamma_p}.
\end{multline*}
In particular, observe that
\begin{equation} \label{eqn:P sub n with lots of zeros}
   \calP_n \big( (0,0,g_0,\bar{g_0},0,0,\dots,0); z, \chi, \hat{z}, \hat{\chi} \big)
    = - g_0 \, S_{\tau^n}(z,\chi,0) + \bar{g_0}^n \, \hat{S}_{\hat{\tau}^n}(\hat{z},\hat{\chi},0).
\end{equation}

This completes the first step of the proof.  We move on to the second step, which involves parametrizing $\Lambda^n$.

\begin{prop} \label{finite det, finally}
Let $(M,0)$ and $(\hat{M},0)$ be formal hypersurfaces of $1$-infinite type which are formally equivalent as above.  Then for every $n \in \N$, there exists a power series $\A_n(z;\Delta, \Lambda) \in \C[\Delta,\Lambda][[z]]^2$ such that the following holds for any $H \in {\F}(M,0;\hat{M},0)$:
\[
   \big( f_n(z), g_n(z) \big) = \A_n \bigg( z; \Delta(H), \big( \lambda^n_{2+\delta^1_K+\delta^1_L \delta^1_T}(H) \big)_{j \in \D(M), j \le n} \bigg).
\]
Moreover, if $M$ and $\hat{M}$ are convergent, then there exists an $\epsilon > 0$ such that
\[
   \bigg\{ z \mapsto \A_n \bigg( z; \Delta(H), \big( \lambda^n_{2+\delta^1_K+\delta^1_L \delta^1_T}(H) \big)_{j \in \D(M), j \le n} \bigg) \bigg\}
   \in \calO_\epsilon(z)^2
\]
for every $n \in \N$ and every $H \in {\F}(M,0;\hat{M},0)$.
\end{prop}

\begin{proof}
We continue with the notation from Proposition \ref{reflection identities}; in particular, we shall continue to let $\gamma$ denote $2 + \delta^1_K + \delta^1_T$.  Observe that Proposition \ref{finite det, finally} follows immediately from Proposition \ref{reflection identities} if it can be shown that for every $n \not \in \D(M)$, there exists a $\C^{\gamma}$-valued polynomial $\omega^n(\Delta, \Lambda)$ such that
\begin{equation} \label{eqn:omega}
   \lambda^n_\gamma(H) = \omega^n \bigg( \Delta(H), \Lambda^{n-1}_{2+\delta^1_K+\delta^1_L \delta^1_T}(H) \bigg) \quad
   \forall \, H \in \F(M,0;\hat{M},0).
\end{equation}
To see this, suppose equation (\ref{eqn:omega}) holds for every $n \not \in \D(M)$.  An easy induction shows that for \emph{every} $n \in \N$, there exists a $\C^\gamma$-valued polynomial $\tilde{\omega}^n(\Delta,\Lambda)$ such that
\[
   \lambda^n_\gamma(H) = \tilde{\omega}^n \bigg( \Delta(H), \big( \lambda^j_{2+\delta^1_K+\delta^1_L \delta^1_T}(H) \big)_{j \in \D, j \le n} \bigg).
\]
Substituting this into the power series for $\B_n$ given by Proposition \ref{reflection identities} completes the proof.

Hence, we must show that a relation of the form given in (\ref{eqn:omega}) holds for each $n \not \in \D(M)$.  To this end, define the power series $\tilde{\Upsilon}^n : (\C^2,0) \to (\C^4,0)$ by $\tilde{\Upsilon}^n_j = \Upsilon^n_j$ for $j \neq 3$, and set
\begin{align*}
   &\tilde{\Upsilon}^n_3(z,\chi) :=
        \delta^1_T \bigg\{ \delta^1_K \bigg( \theta_1^{(L)}(0) \frac{ \theta_\chi(z,\chi) }{ \bar{\theta_L}'(\chi) } \bigg) \\
   &\qquad + \bigg( \frac{ L( \theta_L^{(K)}(0) \theta_{L+1}^{(K)}(0) -\theta_L^{(K+1)}(0) \theta_{L+1}^{(K-1)}(0) ) }{ K \, \theta_L^{(K)}(0)^2 } \bigg) \frac{ \bar{\theta_L}(\chi)}{\bar{\theta_L}'(\chi)} \theta_\chi(z,\chi)  \\
   &\qquad + \bigg( \frac{1+i\,\theta(z,\chi)}{1-i\,\theta(z,\chi)} \bigg)^n \bigg[ \theta_1(z) \big( 1 + \theta(z,\chi)^2 \big) + \bigg( \frac{\theta_{L+1}(z)}{{\theta_L}'(z)} - 2 i \, n \, \frac{\theta_1(z)^2}{{\theta_L}'(z)} \bigg) \theta_z(z,\chi) \bigg]  \\
   &\qquad - \frac{ \theta_{L+1}^{(K-1)}(0) }{ \theta_L^{(K)}(0) } \bigg[ \bar{\theta_1}(\chi) \big( 1 + \theta(z,\chi)^2 \big) + \bigg( \frac{\bar{\theta_{L+1}}(\chi)}{\bar{\theta_L}'(\chi)} + 2 i \, n \, \frac{\bar{\theta_1}(\chi)^2}{\bar{\theta_L}'(z)} \bigg) \theta_\chi(z,\chi) \bigg] \bigg\}
\end{align*}
Observe that $\delta^1_L \, \tilde{\Upsilon}^n_3 = \Upsilon^n_3$.

Reconsider the identity (\ref{eqn:n pair expression}).  If we substitute into it the explicit formulas for $f_n(z)$ and $g_n(z)$ given in Proposition \ref{reflection identities}, as well as the corresponding formulas for $\bar{f_n}(\chi)$ and $\bar{g_n}(\chi)$ given by equation (\ref{eqn:barring the fn and gn's}), we can rewrite this as
\begin{equation} \label{eqn:big linear equation}
   \tilde{\Upsilon}^n(z, \chi)^t \, \kappa^n \big( \Delta(H), \lambda^0_2(H) \big) \, \lambda^n_4(H)
   \equiv W^n \big( z,\chi;\Delta(H), \Lambda^{n-1}_\gamma(H) \big),
\end{equation}
where the superscript $\! \ ^t$ denotes the transpose operation,  $\kappa^n(\Delta,\lambda)$ is the $4 \times 4$ matrix of polynomials defined by
\[
   \kappa^n(\Delta,\lambda^0_2) 
   := \left( \begin{array}{cccc}
       \frac{L}{K} \Delta (b_0^0)^2 & - \frac{n}{K} & - \delta^1_T \frac{L}{K} a_0^2 \, \Delta^2 (b_0^0)^3 & 0 \\
       0 & -\frac{i}{2} & 0 & 0 \\
       0 & 0 & - \delta^1_T \,\Delta (b_0^0)^2 & 0 \\
       0 & 0 & 0 &  \delta^1_K \frac{i}{2}
      \end{array} \right),
\]
(by Proposition \ref{reflection identities 0}, $a_0^2$ is a polynomial in $a_0^1$), and $W^n(z,\chi;\Delta,\Lambda) \in \C[\Delta,\Lambda][[z,\chi]]$.

Denote by $\tilde{\kappa}^n$ the $4 \times 4$ matrix function
\[
   \tilde{\kappa}^n(\Delta,\lambda^0_2)
   := \left( \begin{array}{cccc}
       \frac{K}{L} \Delta (a_0^1)^2 & \frac{2i \, n}{L} \Delta (a_0^1)^2 &  - a_0^2 \, \Delta \, a_0^1 & 0 \\
       0 & 2i & 0 & 0 \\
       0 & 0 & - \delta^1_T \, \Delta (a_0^1)^2 & 0 \\
       0 & 0 & 0 &  - \delta^1_K \, 2i
      \end{array} \right),
\]
Observe that if $a_0^1 \, b_0^0 \neq 0$, then
\[
   \kappa^n \bigg( \frac{1}{a_0^1 \, b_0^0}, \lambda^0_2 \bigg) 
   \cdot \tilde{\kappa}^n \bigg( \frac{1}{a_0^1 \, b_0^0}, \lambda^0_2 \bigg)
   = \left( \begin{array}{cccc}
       1 & 0 & \frac{L \, a_0^2}{K \, a_0^1}(\delta^1_T - 1) & 0 \\
       0 & 1 & 0 & 0 \\
       0 & 0 & \delta^1_T & 0 \\
       0 & 0 & 0 &  \delta^1_K
      \end{array} \right)_j,
\]
For convenience, we shall denote by $\kappa^n_j$ the upper-left $j \times j$ submatrix of $\kappa^n$ for $1 \le j \le 4$; we define $\tilde{\kappa}^n_j$ similarly.  We now complete the proof by examining cases.

\textbf{Case 1:} $K = 1$.  Observe that $L = T = 1$ necessarily, so $\tilde{\Upsilon}^n = \Upsilon^n$ and $\kappa^n_4, \tilde{\kappa}^n_4$ are matrix inverses for all $n \in \N$.  Suppose that $n \not \in \D(M)$, and choose a basis $\{ \upsilon^n_{s_j,t_j} \}_{j=1}^4$ for $\V^n$.  If $\Xi$ is the $4\times4$ matrix whose $j$-th row is $\upsilon^n_{s_j,t_j}$, then it follows that $\Xi$ is invertible.  Now, differentiating (\ref{eqn:big linear equation}) $s_j$ times in $z$, $t_j$ times in $\chi$, and setting $z=\chi=0$ (for $j = 1,2,3,4$), we obtain the $4\times4$ linear system of equations of the form
\[
   \Xi \, \kappa^n_4(\Delta(H),\lambda^2_0(H)) \, \lambda^n_4 = w^n(\Delta(H),\Lambda^{n-1}_4(H)),
\]
Thus, we may take 
\[
   \omega^n(\Delta,\Lambda^{n-1}_4) :=
   \tilde{\kappa}_4^n(\Delta,\lambda^2_0) \, \Xi^{-1} \, w^n(\Delta,\Lambda^{n-1}_4)
\]
to complete the proof.

\textbf{Case 2:} $K > L = 1 = T$.  We have $\tilde{\Upsilon}^n = \Upsilon^n = \big( \Upsilon^n_1,\Upsilon^n_2,\Upsilon^n_3,0)$ and $\kappa^n_3, \tilde{\kappa}^n_3$ are inverses for all $n \in \N$.  Observe too that (\ref{eqn:big linear equation}) reduces to
\begin{multline*}
   \big( \Upsilon^n_1(z, \chi), \Upsilon^n_2(z, \chi), \Upsilon^n_3(z, \chi) \big)^t
      \, \kappa^n_3 \big( \Delta(H), \lambda^0_2(H) \big) \, \lambda^n_3(H)\\
   \equiv W^n \big( z,\chi;\Delta(H), \Lambda^{n-1}_3(H) \big).
\end{multline*}
The proof now follows the exact same lines as in the previous case.

\textbf{Case 3:} $T = 0$.  Since this implies $K > 1$ necessarily, it follows that $\tilde{\Upsilon}^n = \Upsilon^n = \big( \Upsilon^n_1,\Upsilon^n_2,0,0)$ and $\kappa^n_2, \tilde{\kappa}^n_2$ are inverses for all $n \in \N$.  Here, the identity (\ref{eqn:big linear equation}) reduces to
\begin{equation} \label{eqn:the 2 by 2 case}
   \big( \Upsilon^n_1(z, \chi), \Upsilon^n_2(z, \chi) \big)^t \, \kappa^n_2 \big( \Delta(H), \lambda^0_2(H) \big) \, \lambda^n_2(H)
   \equiv W^n \big( z,\chi;\Delta(H), \Lambda^{n-1}_2(H) \big).
\end{equation}
The proof now follows the exact same lines as in the previous two cases.

\textbf{Case 4:} $L > 1 = T$.  Observe that identity (\ref{eqn:big linear equation}) reduces to
\begin{multline} \label{eqn:off a little}
   \big( \Upsilon^n_1(z, \chi), \Upsilon^n_2(z, \chi), \tilde{\Upsilon}^n_3(z, \chi) \big)^t
      \, \kappa^n_3 \big( \Delta(H), \lambda^0_2(H) \big) \, \lambda^n_3(H)\\
   \equiv W^n \big( z,\chi;\Delta(H), \Lambda^{n-1}_3(H) \big).
\end{multline}
We claim that $a_n^0 = \sigma^n(\Delta(H), \Lambda^{n-1}_3(H))$ for every $n \in \N$, where $\sigma^n$ is a polynomial.  Hence, we can write
\[
   \big( f_n(z), g_n(z) \big) = \B_n \big( z;\Delta(H),\Lambda^n_3(H) \big) = \tilde{\B}_n \big( z;\Delta(H),\Lambda^n_2(H) \big);
\]
that is, $f_n(z)$ and $g_n(z)$ are given by expressions of the same form as in Proposition \ref{reflection identities}, \emph{but without the $a_n^0$ term}.  Hence, identity (\ref{eqn:big linear equation}) reduces to identity (\ref{eqn:the 2 by 2 case}), and the proof proceeds as in Case 3.

To prove the claim, we proceed by induction.  For $n = 0$, this is trivial, as $a_0^0 = 0$.  For the inductive step, we consider two cases.

\textbf{Subcase A:} $\theta_{L+1}^{(K-1)}(0) = 0$.  Then equation (\ref{eqn:fn equation formed}) implies
\[
   \bar{a_n^0} = f_n(0) = \bar{a_0^1} \, T^1_n \big( 0;\Delta(H), \Lambda^{n-1}_3(H) \big).
\]
Conjugating this and applying equation (\ref{eqn:barring the fn and gn's}) yields $a_n^0 = \tilde{T}(\Delta(H), \Lambda^{n-1}_3(H))$ for some polynomial $\tilde{T}(\Delta,\Lambda)$.  But by the inductive hypothesis, $\Lambda^{n-1}_3(H)$ is itself a polynomial in $(\Delta(H), \Lambda^{n-1}_2(H))$, so the induction is complete in this case.

\textbf{Subcase B:} $\theta_{L+1}^{(K-1)}(0) \neq 0$.  Differentiating (\ref{eqn:off a little}) $L-1$ times in $\chi$ and setting $\chi = 0$ yields the identity
\[
   \frac{ \big|\theta_{L+1}^{(K-1)}(0) \big|^2}{\big| \theta_L^{(K)}(0) \big|^2} \theta_L(z) \, a_n^0 = W_{\chi^{L-1}} \big( z,0;\Delta(H),\Lambda^{n-1}_3(H) \big).
\]
Differentiating this $K$ times in $z$ and setting $z = 0$ yields $a_n^0 = \tilde{T}(\Delta(H), \Lambda^{n-1}_3(H))$ for some polynomial $\tilde{T}(\Delta,\Lambda)$.  But by the inductive hypothesis, $\Lambda^{n-1}_3(H)$ is itself a polynomial in $(\Delta(H), \Lambda^{n-1}_2(H))$, so the induction is complete in this case.
\end{proof}

This completes the second step.  We move on to the third step, counting the elements of $\D$. 

\begin{prop} \label{big matrix}
Given a fixed set of normal coordinates $Z$ on $M$, the set $\D(M)$ defined by equation (\ref{eqn:D defined}) has at most $2(2+\delta^1_K+\delta^1_L \delta^1_T)$ elements.
\end{prop}

\begin{proof}
Consider the power series $\Upsilon^n(z,\chi)$ defined in equations (\ref{eqn:UP1}) through (\ref{eqn:UP4}); we must prove that for all but $2(2+\delta^1_K+\delta^1_L \delta^1_T)$ integers $n \in \N$, the set $\V^n(M)$ has dimension $2 + \delta^1_K + \delta^1_L \, \delta^1_T$.

Consider the matrix
\[
  \xi(n) := \left(  \begin{array}{cccc}
     \uparrow & \uparrow & \uparrow & \uparrow \\
     \upsilon^n_{2K,2L} & \upsilon^n_{3K,3L} & \upsilon^n_{3K,2L} & \upsilon^n_{2K,3L} \\
     \downarrow & \downarrow & \downarrow & \downarrow 
     \end{array} \right)^t.
\]
Our goal will be to show that for all but at most $2(2+\delta^1_K+\delta^1_L \delta^1_T)$ integers $n \in \N$, the first $2+\delta^1_K+\delta^1_L \delta^1_T$ rows are linearly independent, which implies that $n \not \in \D(M)$.

Using Faa de Bruno's formula, we compute that
\begin{align*}
   (\Upsilon^n_1)_{\chi^{2L}}(z,0) &=  2i \frac{(2L)!}{(L!)^2} K \, \theta_L^{(K)}(z)^2 \, n 
          + \calQ^0 \big( n; \big( \del^\nu \theta(z,0) \big)_{|\nu| < 3L + K + 1} \big)\\
   (\Upsilon^n_1)_{\chi^{3L}}(z,0) &=  -2 \frac{(3L)!}{(L!)^3} K \, \theta_L^{(K)}(z)^3 \, n^2 
          + \calQ^1 \big( n; \big( \del^\nu \theta(z,0) \big)_{|\nu| < 4L + K + 1} \big)\\
   (\Upsilon^n_2)_{\chi^{2L}}(z,0) &=  -2 \frac{(2L)!}{(L!)^2}  \theta_L^{(K)}(z)^2 \, n^2 
          + \calQ^1 \big( n; \big( \del^\nu \theta(z,0) \big)_{|\nu| < 3L + K + 1} \big)\\
   (\Upsilon^n_2)_{\chi^{3L}}(z,0) &=  -  \frac{4i}{3} \frac{(3L)!}{(L!)^3} \theta_L^{(K)}(z)^3 \, n^3 
          + \calQ^2 \big( n; \big( \del^\nu \theta(z,0) \big)_{|\nu| < 4L + K + 1} \big)\\
   (\Upsilon^n_3)_{\chi^{2}}(z,0)  &= \delta^1_L \delta^1_T \big[ - 4 \,  \theta_1^{(K)}(z)^3 \, n^2 
          + \calQ^1 \big( n; \big( \del^\nu \theta(z,0) \big)_{|\nu| < K + 4} \big) \big]\\
   (\Upsilon^n_3)_{\chi^{3}}(z,0)  &= \delta^1_L \delta^1_T \big[ - 16 i \,  \theta_1^{(K)}(z)^4 \, n^3 
          + \calQ^2 \big( n; \big( \del^\nu \theta(z,0) \big)_{|\nu| < K + 5} \big) \big] \\
   (\Upsilon^n_4)_{\chi^{2}}(z,0)  &= \delta^1_K \big[ 
            \calQ^0 \big( n; \big( \del^\nu \theta(z,0) \big)_{|\nu| <  5} \big) \big] \\
   (\Upsilon^n_4)_{\chi^{3}}(z,0)  &= \delta^1_K \big[ 12 \, \theta_1(z)^2 \, n^2 
          + \calQ^1 \big( n; \big( \del^\nu \theta(z,0) \big)_{|\nu| < 5} \big) \big].
\end{align*}  
Setting $\alpha := \theta_L^{(K)}(0)$ it follows, we may write $\xi(n) = C_1(n) + C_2(n)$, with
\[
   C_1(n) = \left( \begin{array}{cccc}
    \frac{2 i\,K (2L)!(2K)! \alpha^2}{(L!K!)^2}  n &  \frac{-2(2L)!(2K)! \alpha^2 }{(L!K!)^2} n^2 & 0 & 0 \\
    \frac{-2 K(3L)!(3K)! \alpha^3 }{(L!K!)^3} n^2 & \frac{-4i(3L)!(3K)! \alpha^3}{3(L!K!)^3}  n^3 & 0 & 0 \\
    0 & 0 & \delta^1_L \delta^1_T \frac{- 4 (3K)! \alpha^3}{(K!)^3} n^2 & 0 \\
    0 & 0 & 0 & \delta^1_K \, 72 \, \alpha^2 n^2 
   \end{array} \right)
\]
and $C_2(n)$ of the form
\[
  \left( \begin{array}{cccc}
    \calQ^0(n; j^{3L+3K+1}_0 \theta) & \calQ^1(n; j^{3L+3K+1}_0 \theta) 
       & \delta^1_L \, \delta^1_T \, \calQ^1(n; j^{3K+4}_0 \theta) & \delta^1_K \, \calQ^0(n; j^{7}_0 \theta) \\
    \calQ^1(n; j^{4L+4K+1}_0 \theta) & \calQ^2(n; j^{4L+4K+1}_0 \theta) 
       & \delta^1_L \, \delta^1_T \, \calQ^2(n; j^{4K+5}_0 \theta) & \delta^1_K \, \calQ^2(n; j^{9}_0 \theta) \\
    \calQ^1(n; j^{3L+4K+1}_0 \theta) & \calQ^2(n; j^{3L+4K+1}_0 \theta) 
       & \delta^1_L \, \delta^1_T \, \calQ^1(n; j^{4K+4}_0 \theta) & \delta^1_K \, \calQ^0(n; j^{8}_0 \theta) \\
    \calQ^1(n; j^{4L+3K+1}_0 \theta) & \calQ^2(n; j^{4L+3K+1}_0 \theta) 
       & \delta^1_L \, \delta^1_T \, \calQ^2(n; j^{3K+5}_0 \theta) & \delta^1_K \, \calQ^2(n; j^{8}_0 \theta) 
   \end{array} \right).
\]
We shall denote by $\xi_j(n)$ the upper-left $j \times j$ submatrix of $\xi(n)$ for $j = 1,2,3,4$.  We complete the proof by examining cases.

\textbf{Case 1:} $K = 1$.  In this case $L = T = 1$ as well, whence $2 + \delta^1_K + \delta^1_L\delta^1_T = 4$.  By examining the matrix $\xi_4(n)$, and in particular  the term of highest order in $n$ in each of its entries, we find that
\[
  \det \xi_4(n) = 110592 \, \alpha^{10} \, n^8 + \calQ^7 \big( n;j^9_0 \theta \big).
\]
Since $\alpha \neq 0$, this is a nonzero, eighth degree polynomial in $n$, and hence has at most eight distinct zeros (in the complex plane).  If $\det \xi_4(n_0) \neq 0$, then the four rows of $\xi(n_0)$ are linearly independent,  which completes the claim.

\textbf{Case 2:} $K > L = T = 1$.  In this case, we have  $2 + \delta^1_K + \delta^1_L\delta^1_T = 3$.  By examining the highest order terms in $n$ as above, we find that
\[
   \det \xi_3(n) =  64 \, K \frac{(2K)!(3K)!^2}{(K!)^8} \alpha^8 \, n^6 + \calQ^5 \big (n;j^{4K+5}_0 \theta \big).
\]
Arguing as above implies that for all but (at most) six integers $n$, the matrix $\xi_3(n)$ is invertible, whence the first three rows of $\xi(n)$ are linearly independent. This completes the claim.

\textbf{Case 3:} $L > 1$ \textbf{or} $T = 0$.  Since either of these conditions necessarily implies $K > 1$, we conclude that $2 + \delta^1_K + \delta^1_L\delta^1_T = 2$. Since
\[
   \det \xi_2(n) = - \frac{4}{3} K \frac{(2L)!(3L)!(2K)!(3K)!}{(L! K!)^5} \alpha^5 \, n^4 + \calQ^3 \big(n;j^{4L+4K+1}_0 \theta \big),
\]
the proof is complete by arguments similar to the previous case.
\end{proof}

It is worthwhile to note that while $\D(M)$ is always finite, it is also never empty.  Indeed, $0 \in \D(M)$ for any $1$-infinite type hypersurface $M$, since it is easy to check that $\Upsilon^0_2(z,\chi) \equiv 0$.

This completes the third step of the proof.  We complete the proof by showing that $\D(M)$ is independent of the choice of normal coordinates used to define it.  In fact, we prove the following, which completes the proof of Theorem \ref{linear independence}.

\begin{prop} \label{invariance}
Suppose that $M$, $Z = (z,w)$, $\Upsilon^n$, and $\V^n = \V^n(M)$ are as above. Let $(\hat{M},0)$ be formally equivalent to $(M,0)$, with corresponding power series $\hat{\Upsilon}^n$ and subspaces $\hat{\V}^n = \V^n(\hat{M})$ defined using the normal coordinates $\hat{Z}=(\hat{z},\hat{w})$.  Then for every $n \in \N$, the dimensions of $\V^n$ and $\hat{\V}^n$ are equal.  In particular, the dimension of subspace $\V^n(M) \subset \C^4$ is independent of the choice of normal coordinates used to define it.
\end{prop}

\begin{proof}
Let $H(z,w) = \big( f(z,w), w \, g(z,w) \big)$ be a formal equivalence between $M$ and $\hat{M}$.  Consider the formal power series 
\[
   (z,\chi) \mapsto \hat{\Upsilon}^n \big( f_0(z), \bar{f_0}(\chi) \big)
   \in \C [[ z, \chi ]]^4,
\]
which may be viewed as the power series $\hat{\Upsilon}^n$ given in the $Z$ coordinates.  Using Faa de Bruno's formula and the fact that $f_0 : (\C,0) \to (\C,0)$ is a formal change of coordinates, it is straightforward to verify that
\[
   \Span_\C \setof{ \hat{\upsilon}^n_{s,t} := \dd{^{s+t}}{z^s \del \chi^t} \big\{ \hat{\Upsilon}^n \big( f_0(z), \bar{f_0}(\chi) \big) \big\} \bigg|_{ \substack{z = 0 \\ \chi = 0} } }{ s, t \in \N }
   = \hat{\V}^n.
\]

From (\ref{eqn:f0 in phi}) we derive
\[
   \hat{\theta}_{\hat{z}}  \big( f_0(z), \bar{f_0}(\chi) \big) = \frac{\theta_z(z,\chi)}{{f_0}'(z)}, \quad
   \hat{\theta}_{\hat{\chi}}  \big( f_0(z), \bar{f_0}(\chi) \big) = \frac{\theta_\chi(z,\chi)}{\bar{f_0}'(\chi)}
\]
whereas repeated differentiation of this in $\chi$ yields
\[
   \hat{p}_{L+1} \big( f_0(z) \big) = \frac{1}{2 (a_0^1)^{L+2}} \big( 2 a_0^1 \, p_{L+1}(z) - (L+1)L \, a_0^2 \, p_L(z) \big).
\]
From this and identity (\ref{eqn:f0 fuller}), it follows by an elementary (albeit involved) calculation that
\[
\begin{array}{r@{\,\, = \,\,}c@{\,}l}
   \hat{\Upsilon}^n_1 \big( f_0(z), \bar{f_0}(\chi) \big) & &
        \Upsilon^n_1(z,\chi)\\
   \hat{\Upsilon}^n_2 \big( f_0(z), \bar{f_0}(\chi) \big) & & 
        \Upsilon^n_2(z,\chi)\\
   \hat{\Upsilon}^n_3 \big( f_0(z), \bar{f_0}(\chi) \big) & \frac{\displaystyle 1}{\displaystyle a_0^1} &
        \Upsilon^n_3(z,\chi) 
        + \frac{\displaystyle \delta^1_T  \, a_0^2}{\displaystyle K(a_0^1)^2} \Upsilon^n_1(z,\chi)\\
   \hat{\Upsilon}^n_4 \big( f_0(z), \bar{f_0}(\chi) \big) & a_0^1 & 
        \Upsilon^n_4(z,\chi)
\end{array}
\]

Now, suppose that $\{ \hat{\upsilon}^n_{s_j,t_j} \}_{j=1}^{\ell_0}$ is any collection of vectors in $\hat{\V}^n$; consider the corresponding vectors $\upsilon^n_{s_j,t_j} \in \V^n$.  Observe that if $\hat{\Xi}, \Xi$ denote the $4 \times \ell_0$ matrices whose columns are, respectively, the $\hat{\upsilon}^n_{s_j,t_j}, \upsilon^n_{s_j,t_j}$, then in view of the above identities, these matrices necessarily have the same rank.  In particular, the columns of $\hat{\Xi}$ are linearly independent if and only if the columns of $\Xi$ are.  From this it follows that $\hat{\V}^n$ and $\V^n$ have the same dimension.
\end{proof}

\subsection{The main results}

We use Theorem \ref{linear independence} to prove the main theorems stated at the end of Section \ref{sec:Preliminaries and Basic Definitions}.  We begin with Theorem \ref{Parametrization}.  

\begin{proof}  
Let $M$ be a formal real hypersurface of $1$-infinite type at $0$.  Observe that the result of Theorem \ref{Parametrization} is independent of the choice of coordinates $Z$, so without loss of generality let us take $Z = (z,w)$ to be normal coordinates for $M$, so that $M$ is given by equation (\ref{eqn:I}).  Let $\D = \D(M)$ be as in Theorem \ref{linear independence}, and set $k := 2 + \max \D$, which exists since $\D$ is a finite set.

To prove this $k$ is sufficient, suppose $\hat{M}$ is a formally equivalent formal real hypersurface.  Define the corresponding $\A^n$ as in Theorem \ref{linear independence}.  Fix a formal equivalence $H \in \F(M,0;\hat{M},0)$.  Conjugating the formula for $(f_n,g_n)$ implies that
\[
   \big( \bar{f_n}(\chi), \bar{g_n}(\chi) \big) =
   \bar{\A}^n \bigg( \frac{1}{ \bar{a_0^1} \, \bar{b_0^0} }, \big(\bar{a_j^0}, \bar{b_j^0}, \bar{a_j^1}, \bar{b_j^1} \big)_{j \in \D} \bigg),
\]
whence
\[
   (a_n^0,b_n^0,a_n^1,b_n^1) = A_n \bigg( \frac{1}{ \bar{a_0^1} \, \bar{b_0^0} }, \big(\bar{a_j^0}, \bar{b_j^0}, \bar{a_j^1}, \bar{b_j^1} \big)_{ j \in \D} \bigg), \quad n=\N
\]
with $A_n \in \C[Delta,\Lambda]^4$.  Substituting this into $\A_n$ --- and recalling that 
\[
   \Delta(H) = \frac{1}{a_0^1 \, b_0^0} = \frac{\bar{a_0^1}}{ \mu^2 \bar{b_0^0}}, 
\]
where $\mu$ is defined by (\ref{eqn:a01 constant norm}) --- we can write
\[
   \big( f_n(z), g_n(z) \big) = \Gamma^n \bigg( z; \frac{1}{ \bar{a_0^1} \, \bar{b_0^0} }, \big(\bar{a_j^0}, \bar{b_j^0}, \bar{a_j^1}, \bar{b_j^1} \big)_{j \in \D} \bigg),
\]
with $\Gamma^n(z;\Delta,\Lambda) \in \C[\Delta,\Lambda][[z]]^2$.  Let us write
\[
   \Gamma^n_{z^j} \bigg( 0; \frac{1}{ \bar{a_0^1} \, \bar{b_0^0} }, \big(\bar{a_j^0}, \bar{b_j^0}, \bar{a_j^1}, \bar{b_j^1} \big)_{j \in \D} \bigg)
   =: \frac{ c^n_j \bigg( \big(\bar{a_j^0}, \bar{b_j^0}, \bar{a_j^1}, \bar{b_j^1} \big)_{j \in \D} \bigg) }{ \big( \bar{a_0^1} \, \bar{b_0^0} \big)^{\ell^n_j} },
\]
with $\ell^n_j \in \N$ and $c_j^n$ a $\C^2$-valued polynomial.

Now, observe that
\[
   \dd{^{\ell+j} H}{z^\ell \, \del w^j}(0,0) = \bigg( \, \bar{a_j^\ell}, j \, \bar{b_{j-1}^{\ell}} \, \bigg).
\]
In particular, observe that $\bar{a_j^0}$ is a term in (the coordinates of) $j^k_0(H)$, $a_j^1$ and $b_j^0$ are terms in $j^{k+1}_0(H)$, and $b_j^1$ is a term in $j^{j+2}_0(H)$.  Hence, $c_n^j$ is a polynomial in $j^{2 + \max \D}_0 (H) = j_0^k(H)$ and
\[
   0 \neq \bar{a_0^1} \, \bar{b_0^0} = \det \bigg( \dd{H}{Z}(0,0) \bigg) =: q \big( j_0^k (H) \big)
\]
so the proof is complete in view of equation (\ref{eqn:fn and gn expanded}).
\end{proof}

Observe by inspecting Propositions \ref{reflection identities} through \ref{invariance}, we can actually replace the $k$ given in the proof by $k := 1 + \delta^1_K + \max \D$ to get a better bound in the $K > 1$ case, and if $\D = \{ 0 \}$, then we may take $k = 1$ since $b_0^1 = 0$ by Corollary \ref{f0 equation found}.

We now use this result to prove Theorem \ref{Finite Determination}.

\begin{proof}
Let $M, k$ be as in Theorem \ref{Parametrization}.  Suppose that $\hat{M}$ is formally equivalent to $M$, and let $\Psi$ be the formal power series from Theorem \ref{Parametrization}.  If $H^1, H^2 : (M,0) \to (\hat{M},0)$ are two formal equivalences which satisfy
\[
   \dd{^{|\alpha|} H^1}{Z^\alpha}(0) = \dd{^{|\alpha|} H^2}{Z^\alpha}(0)
   \quad \forall \, |\alpha| \le k,
\]
then it follows that $j^k_0 (H^1) = j^k_0 (H^2)$.  If we call this common jet $\Lambda_0$, then it follows from Theorem \ref{Parametrization} that
\[
   H^1(Z) \equiv \Psi(Z; {\Lambda_0}) \equiv H^2(Z),
\]
as desired.
\end{proof}

We now tackle the two applications of Theorem \ref{Parametrization} mentioned in Section \ref{sec:Preliminaries and Basic Definitions}.  First we prove Theorem \ref{Lie Subgroup}.

\begin{proof}
Let $M,k$ be as in Theorem \ref{Parametrization}, and let $\Psi$ be the formal power series defined in accord with that theorem with $\hat{M} = M$.  That the mapping $j^k_0 : \Aut(M,0) \to J^k_0(\C^2,\C^2)_{0,0}$ is injective follows from Theorem \ref{Finite Determination}. Observe that $\Lambda_0 \in J^k(\C^2,\C^2)_{0,0}$ is in the image of $j^k_0$ if and only if $q(\Lambda_0) \neq 0$ --- so that $\Lambda_0 \in G^k(\C^2)_0$) --- and
\begin{gather}
  \label{eqn:injective}
   \Lambda_0 = j^k_0 \big( \Psi(\cdot,\Lambda_0) \big) \\
  \label{eqn:H(M) to M'}
   \rho \big( \Psi(Z,\Lambda_0), \bar{\Psi}(\zeta,\bar{\Lambda_0}) \big) = a(Z,\zeta) \rho(Z,\zeta) 
\end{gather}
for some multiplicative unit $a(Z,\zeta) \in \C[[Z,\zeta]]$, where $\rho$ is a defining power series for $M$.  In view of equation (\ref{eqn:form of Psi II}),  (\ref{eqn:injective}) is a finite set of polynomial equations in $\Lambda_0$, whereas (\ref{eqn:H(M) to M'}) is a (possibly countably infinite) set of polynomial equations in $(\Lambda_0,\bar{\Lambda_0})$.  Hence, the image of the mapping $j_0^k$ is a locally closed subgroup of the Lie group $G^k(\C^2)_0$, and so is a Lie subgroup.
\end{proof}

And as a corollary, we have Theorem \ref{Varieties}.

\begin{proof}
Let $M, k$ be as in Theorem \ref{Parametrization}, and let $(\hat{M},0)$ be formally equivalent to $(M,0)$.  Injectivity of the jet map again follows from Theorem \ref{Finite Determination}.  Now, fix a formal equivalence $H_0 : (M,0) \to (\hat{M},0)$; then any other formal equivalence is of the form $H := H_0 \circ A$, where $A \in \Aut(M,0)$.  In particular,
\begin{align*}
   j^k_0 \big( {\F}(M,0;\hat{M},0) \big)
    &= \setof{ j^k_0(H_0 \circ A) }{ A \in \Aut(M,0) } \\
    &= \setof{ j^k_0(H_0) \cdot j^k_0(A) }{ A \in \Aut(M,0) } \\
    &= j^k_0(H_0) \cdot j^k_0 \big( \Aut(M,0) \big).
\end{align*}
Hence, the image of ${\F}(M,0;\hat{M},0)$ is merely a coset of the algebraic Lie subgroup $j^k_0 \big( \Aut(M,0) \big)$ in the Lie group $G^k(\C^2)_0$, and so is itself a real-algebraic submanifold of $G^k(\C^2)_0$. 
\end{proof}

\end{document}